\documentclass[reqno,10pt,draft,a4paper]{amsart}
\usepackage{amsmath}
\usepackage{amssymb}
\usepackage{latexsym}
\usepackage{epsf}

\usepackage[cmtip,matrix,arrow]{xy}
\UseComputerModernTips
\CompileMatrices

\DeclareMathOperator{\Hom}{Hom}
\DeclareMathOperator{\Ext}{Ext}

\renewcommand{\ge}{\geqslant}
\renewcommand{\le}{\leqslant}

\newcommand{\N}{\mathbb{N}}
\newcommand{\R}{\mathbb{R}}

\newcommand{\Z}{\mathbb{Z}}

\newcommand{\frob}{^{\mathrm F}}

\newcommand{\rarr}{\rightarrow}

\DeclareMathOperator{\pr}{pr}

\newcommand{\lrimpl}{\Leftrightarrow}
\DeclareMathOperator{\hd}{hd}
\DeclareMathOperator{\soc}{soc}
\newcommand{\Mod}{\mathrm{mod}}
\newcommand{\MMod}{\mathrm{Mod}}
\DeclareMathOperator{\ch}{ch}

\DeclareMathOperator{\Ind}{Ind}

\newcommand{\ses}{short exact sequence }
\newcommand{\GL}{\mathrm{GL}}
\newcommand{\SL}{\mathrm{SL}}
\newcommand{\St}{\mathrm{St}}
\newcommand{\half}{\frac{1}{2}}

\newcommand{\crho}{\rho\check{\ }\,}
\newcommand{\alphac}{\alpha\check{\ }\,}

\newcommand{\betac}{\beta\check{\ }\,}
\newcommand{\bnabla}{\overline{\nabla}}

\begin{document}
\theoremstyle{plain}
\numberwithin{subsection}{section}
\newtheorem{thm}{Theorem}[section]
\newtheorem{propn}[thm]{Proposition}
\newtheorem{cor}[thm]{Corollary}
\newtheorem{clm}[thm]{Claim}
\newtheorem{lem}[thm]{Lemma}
\newtheorem{conj}[thm]{Conjecture}
\theoremstyle{definition}
\newtheorem{defn}[thm]{Definition}
\newtheorem{rem}[thm]{Remark}

\title{Good $l$-filtrations for $q-\GL_3(k)$}
\author{Alison E. Parker}
\address{School of Mathematics and Statistics F07, University of Sydney,
  NSW 2006, Australia}
\email{alisonp@maths.usyd.edu.au}




\begin{abstract}
Let $k$ be an algebraically closed field of characteristic $p$,
possibly zero, and  $G=q$-$\GL_3(k)$, the
quantum group of three by three matrices as defined by
Dipper and Donkin. We may also take $G$ to be
$\GL_3(k)$. 
We first determine the extensions between simple $G$-modules 
for both $G$ and $G_1$, the first
Frobneius kernel of $G$. 
We then determine the submodule structure of certain
induced modules, $\hat{Z}(\lambda)$, for the infinitesimal group $G_1B$.
We induce this structure to $G$ to obtain a good $l$-filtration
of certain induced modules, $\nabla(\lambda)$, for $G$. 
We also determine the homomorphisms between induced modules for $G$.
\end{abstract}
\maketitle

\section*{Introduction}
Let $k$ be an algebraically closed field of characteristic $p$,
possibly zero.
In this paper we study the module category for $G=q$-$\GL_3(k)$, the
quantum group of three by three matrices. We use the quantisation of
Dipper and Donkin \cite{dipdonk}. We may also take $G$ to be
$\GL_3(k)$, that is the classical group scheme of three by three
invertible matrices.

We want to determine explicitly the structure of two types of
modules. First we determine the submodule structure of certain
induced modules, $\hat{Z}(\lambda)$, for the infinitesimal group $G_1B$.
We then induce this structure to $G$ to obtain a good $l$-filtration
of certain induced modules, $\nabla(\lambda)$, for $G$. 
We also determine the homomorphisms between induced modules for $G$.

This paper generalises several classical results
including the extensions between simple modules for $\SL_3(k)$,
\cite{yehia}, the submodule structure of the $\hat{Z}(\lambda)$'s for
$\SL_3(k)$, \cite{irv}, some results about translations, \cite{jantz2},
good $p$-filtrations of the induced modules $\nabla(\lambda)$
for $\SL_3(k)$, \cite{parker1}, and the homorphisms between induced
modules for $\SL_3(k)$, \cite{coxpar}.
It also clears up some
confusion regarding the validity of results of Irving \cite{irv}
and Parker \cite{parker1} for small primes.
A large part of this paper produces a quantum version of many results
of the PhD thesis of Yehia, \cite{yehia}. We have reproduced some of his
arguments, only applied to the quantum case, as this reference is not
that accessible. 

\section{Notation}\label{sect:notn}
We first review the basic concepts and most of the notation that we will be using.
A very brief introduction to the theory of quantum groups and how it
relates to linear algebraic groups may
be found in~\cite[chapter 0]{donkbk}.
Some of the cohomological theory of quantum groups and their $q$-Schur
algebras appears in~\cite{donkquant}. We will also refer to
\cite{andpolwen} for many of the basic properties of quantum groups.

Throughout this paper $k$ will be an algebraically closed field of
characteristic $p$ which may be zero.

First take $G$ to be $\GL_3(k)$.
We take 
$l$ to be $p$ which we assume for this particular case to be non-zero. 
We let $T$ be the diagonal matrices in $G$ and $B$, a Borel
subgroup, be the lower triangular matrices. 
We will write $\MMod(G)$ for the category of 
dimensional rational $G$-modules and 
$\Mod(G)$ for the category of finite
dimensional rational $G$-modules. 
We let $D$ be the one-dimensional determinant module for $G$.

Now take $G$ to be $q$-$\GL_3(k)$ the quantum group of Dipper and
Donkin, as defined in \cite{donkbk}.
We write $\MMod(G)$ for the category of
right comodules of $k[G]$, the Hopf algebra of $G$
and $\Mod(G)$ for the category of finite dimensional right comodules
of $k[G]$.
If $q$ is not a root of unity then $\Mod(G)$ is semi-simple.
We will thus consider the case where $q$ is a primitive $l$th root of unity with 
$l \ge 2$.
We take $T$, and $B$ as defined in
\cite{donkquant}.
We let $D$ be the one-dimensional module for $G$, where $G$ acts by
the quantum determinant as defined in \cite{donkquant}.

We now consider both cases together.

Let $X(T)=X \cong \Z^3$ be the weight lattice for $G$ with $\Z$-basis
$\{ e_1=(1,0,0), e_2=(0,1,0), e_3=(0,0,1) \}$. 
Every module in $\Mod(G)$ is semi-simple as a $T$-module and we define
the formal character $\ch (V) \in \Z X$ of $V$ to be the character of
$V$ restricted to $T$. We use $e(\lambda)$ with $\lambda \in X$ as a
basis for $\Z X$, so to distinguish characters from the structure of
the weight lattice as a $\Z$ vector space. We thus have $e(\lambda)
e(\mu) = e(\lambda +\mu)$ in $\Z X$.

We set $R= \{ e_i -e_j \mid i \ne j \}$ to be the roots of $G$. 
For each $\alpha \in R$ we take 
$\alphac=\alpha \in X$ to be  the coroot of $\alpha$. 
(Here we have identified the weight space with the dual weight space,
as we are only considering $\GL_3$, the two are isomorphic.)
Let  $R^+=\{  e_i -e_j \mid i < j \}$ be the positive roots, 
(chosen so that $B$ is the negative
 Borel) and let $S = \{ e_i -e_{i+1} \}$ be the set of simple roots. 
Set $\rho = \half \sum_{\alpha \in R^+} \alpha = (1,0,-1)$.

We have a partial order on $X$ defined by 
$\mu \le \lambda \lrimpl \lambda -\mu \in \N S$.
We also have a bilinear form $\langle -,- \rangle :X  \times X
\rarr \Z$ with $\langle e_i, e_j\check{\ } \rangle = \delta_{ij}$ (Kronecker
delta).
A weight $\lambda$ is \emph{dominant} if
$\langle \lambda, \alphac \rangle \ge 0$ for all $\alpha \in S$ and
we let $X^+$ be the set of dominant weights.
In this case $X^+ = \{ (a,b,c) \mid a \ge b \ge c \}$.

Take $\lambda \in X^+$ and let $k_\lambda$ be the one-dimensional module
for $B$ which has weight $\lambda$. We define the induced
module, $\nabla(\lambda)= \Ind_B^G(k_\lambda)$. 
This module has formal character given by Weyl's character formula and has
simple socle $L(\lambda)$, 
the irreducible $G$-module of highest weight
$\lambda$.  These completely exhaust the simple modules in $\Mod(G)$.
We will denote the socle of a module $M$ by $\soc(M)$.


We return to considering the weight lattice $X$ for $G$.
We consider the affine reflections
$s_{\alpha,ml}$ for $\alpha$ a positive
root and $m\in \Z$ which act on $X$ as
$s_{\alpha,ml}(\lambda)=\lambda -(\langle\lambda,\alphac\rangle -ml )\alpha$.
These generate the affine Weyl group $W_l$. 
We let $W$ be the Weyl
group of $G$ which is generated by $s_{(1,-1,0), 0}$ and 
$s_{(0,1,-1),0}$. 
We mostly use the dot action of $W_l$ on $X$ which is 
the usual action of $W_l$, with the origin
shifted to $-\rho$. So we have $w \cdot \lambda = w(\lambda+\rho)-\rho$.
The reason for this is the following, sometimes known as the linkage
principle.

\begin
{propn}[{\cite[corollary 8.2] {andpolwen}}]
Let $V \in \Mod(G)$ and $V$ be indecomposable. If $L(\mu)$ and
$L(\lambda)$ are composition factors of $V$ then $\mu \in W_l \cdot
\lambda$.
\end{propn}

We now define the quantum version of translation functors. These are defined
in~\cite[section 8]{andpolwen}. 
For any $G$-module $V$ and any $\mu\in X$, set
$\pr_\mu V$ equal to the sum of submodules of $V$ such that all the
composition factors have highest weight in $W_p \cdot \mu$. Then
$\pr_\mu V$ is the largest submodule of $V$ with this property.

\begin{defn}
Suppose $ \lambda$, $\mu\in \bar{C}$. There is a unique $\nu_1 \in
X^+\cap W(\mu-\lambda)$. We define the \emph{translation functor}
$T_\lambda ^\mu$ from $\lambda$ to $\mu$ via
$$T_\lambda ^\mu V = \pr _\mu (L(\nu_1)\otimes \pr_\lambda V)$$
for any $G$-module $V$. It is a functor from $\Mod(G)$ to
itself.
\end{defn}
These functors have similar properties to the classical ones, as
remarked in \cite[section 8]{andpolwen}. 


A \emph{facet} for $W_l$ is a non-empty set of the form
\begin{equation*}
\begin{split}
F= \{ \lambda \in X\otimes_{\Z}\R\ \mid \ 
&\langle \lambda+\rho, \alphac
\rangle = n_\alpha l\quad \forall\, \alpha \in R^+_0(F),
\\
&(n_\alpha -1)l < \langle \lambda +\rho, \alphac \rangle < n_\alpha l
\quad \forall\, \alpha \in R_1^+(F)\} 
\end{split}
\end{equation*}
for suitable $n_\alpha \in \Z$ and for a disjoint decomposition
$R^+=R_0^+(F) \cup R_1^+(F)$.

The \emph{closure} $\bar{F}$ of a facet $F$ is similar but with the inequalities
replaced with equalities.
The \emph{upper  closure} $\hat{F}$ of a facet $F$ is defined as 
\begin{equation*}
\begin{split}
\hat{F}= \{ \lambda \in X\otimes_{\Z}\R\ \mid \ 
&\langle \lambda_\rho, \alphac
\rangle = n_\alpha l\quad \forall\, \alpha \in R^+_0(F),
\\
&(n_\alpha -1)l < \langle \lambda +\rho, \alphac \rangle \le n_\alpha l
\quad \forall\, \alpha \in R_1^+(F)\} 
\end{split}
\end{equation*}

A facet $F$ is an \emph{alcove} if $R^+_0(F)=\emptyset$, (or
equivalently $F$ is open in $X\otimes_{\Z} \R$).
If $F$ is an alcove for $W_l$ then its closure $\bar{F}\cap X$ is a
fundamental domain for $W_l$ operating on $X$. The group $W_l$
permutes the alcoves simply transitively.
We set 
$C= \{ \lambda \in X \otimes _{\Z} \R \ \mid\  0< \langle \lambda +\rho,
\alphac \rangle < l \quad \forall\, \alpha \in R^+ \}$ 
 and call $C$ the \emph{fundamental alcove}.
We have
$C \cap X \ne \emptyset$ if and only if  $l \ge 3$, the Coxeter number
of $G$.

A facet $F$ is a \emph{wall} if 
there exists a unique $\beta\in R^+$ with $\langle \lambda
+\rho, \betac \rangle =ml$ for some $m \in \Z$ and for all $\lambda \in F$.
%

The category $\MMod(G)$ has enough injectives and so we may
define $\Ext_G^*(-,-)$ as usual by using injective
resolutions (see \cite{benson}, section 2.4 and 2.5).

We let ${\mathrm{F}}$ be the Frobenius morphism from $G \to \GL_3(k)$,
and denote by $M\frob$ 
the Frobenius twist of a
module for $\GL_3(k)$. We will sometimes distinguish modules for
$\GL_3(k)$ and $G$ by a bar, $\bar{\ }$.
We set $X_1$ to be the $l$-restricted weights. Thus
$X_1 = \{ (\lambda_1, \lambda_2, \lambda_3) \mid 
0 \le \lambda_1-\lambda_2 < l \mbox{ and }
0 \le \lambda_2-\lambda_3 < l  \}$.
We let $G_1$ be the kernel of $F$ as a group scheme, (it has defining
ideal generated by $c_ij^l - \delta_{ij}$ where the $c_{ij}$ are the
coordinate functions generating
the Hopf algebra $k[G]$ and  $\delta_{ij}$ is the Kronecker delta).
 
We define $\lambda'$ and $\lambda''$ for $\lambda \in X^+$,  
$\lambda = l \lambda''+ \lambda' $  with $\lambda'' \in X^+$
and $\lambda' \in X_1$.
We will use Steinberg's tensor product theorem:
$L(\lambda) \cong \bar{L}(\lambda'')\frob \otimes
L(\lambda')$,
where $\lambda \in X^+$.
We define $\nabla_l(\lambda) = \bnabla(\lambda_1)\frob \otimes
L(\lambda_0)$.

We let $\hat{Z}(\lambda) = \Ind_{B}^{G_1B} k_\lambda$ and 
$\hat{L}(\lambda)$ be the simple module for $G_1B$ of highest weight
$\lambda$. (Note: this is the $\hat{Z}'(\lambda)$ of \cite{jantz2}, we
have dropped the primes, and so our $\hat{Z}(\lambda)$ is
not to be confused with the $\hat{Z}(\lambda)$ of \cite{jantz2}.
The subgroup $G_1B$ has defining ideal generated by $c^l_{ij}$ with $i<j$. Our
$\hat{Z}(\lambda)$ upon restriction to $G_1T$, the subgroup with
defining ideal generated by $c^l_{ij}$ with $i \ne j$, 
is the $\hat{\nabla}_1(\lambda)$ of \cite{donkbk}
and our ``$G_1T$'' is the Janzten subgroup $\hat{G}_1$ of
\cite{donkbk}. This reference doesn't consider the case with
$G_1B$. But many properties for $G_1B$ can be deduced from the
properties for $G_1T$.)
We have $\hat{L}(\lambda) \cong L(\lambda') \otimes k_{l\lambda''}$. 
We will often use a hat $\hat{\ }$ to distinguish modules for $G_1B$
from those for $G$.
Note that we have $\nabla_l(\lambda) \cong \Ind_{G_1B}^G(\hat{L}(\lambda))$. 
We also note that the $\nabla_l(\lambda)$ are indecomposable with simple socle
$L(\lambda)$.

We denote the composition multiplicity of a simple module $L$ in a
module $M$ by $[M:L]$.

Suppose a $G$-module $M$ has a filtration:
$$ 0 = M_0 
\subseteq M_1 \subseteq \cdots \subseteq M_{m-1} \subseteq M_m, $$
with quotients $Q_i= M_i/M_{i-1}$. This will be depicted graphically
as
$$
\xymatrix@R=15pt{
*=0{\bullet} \ar@{-}[d]^{\textstyle{
{Q_m}}} \\
*=0{\bullet} \ar@{-}[d]^{\textstyle{
{Q_{m-1}}}} \\
*=0{\bullet} \ar@{.}[dd]^{\textstyle{
{}}} \\
*=0{}
 \\
*=0{\bullet} \ar@{-}[d]^{\textstyle{
{Q_2}}}\\
*=0{\bullet} \ar@{-}[d]^{\textstyle{
{Q_1}}}\\
*=0{\bullet} 
}
$$
We will also draw pictures like so
$$
\xymatrix@R=15pt{
{Q_m} \ar@{-}[d]\ar@{-}[dr]
&
& 
{Q_{m-1}} \ar@{-}[dl]\\
{Q_{m-2}} 
& 
{Q_{m-3}} 
&
\\
{\vdots} &{\vdots} & {\vdots}
\\
{Q_2}\ar@{-}[dr]
&&
{Q_3}\ar@{-}[dl]\\
&
{Q_1} &\\
}
$$
when we have more information about the extensions appearing between
the $Q_i$ in the module $M$. So the above picture represents a module
with an indecomposable submodule with $Q_1$ and $Q_2$ as factors, etc.

If every quotient $Q_i$ is isomorphic to $\nabla(\mu_i)$ for
some $\mu_i \in X^+$ then we say that $M$ has a \emph{good filtration}.
If every quotient $Q_i$ is isomorphic to $\nabla_l(\mu_i)$ for
some $\mu_i \in X^+$ then we say that $M$ has a \emph{good 
$l$-filtration}. We will often abbreviate this to just $l$-filtration.
If every quotient $Q_i$ is isomorphic to dual induced modules 
$\nabla(\mu_i)^*$ for
some $\mu_i \in X^+$ then we say that $M$ has a \emph{Weyl filtration}.

Good filtration multiplicities and Weyl filtration multiplicities, 
like composition multiplicities are
well defined. It is conjectural that the same holds for good
$l$-filtration multiplicities. They are if a
conjecture of Donkin holds --- this is the subject of
\cite{anderpfil}.

We say a module is a \emph{tilting module} if it has both a good
filtration and a Weyl filtration. For each $\lambda \in X^+$ there is
a unique indecomposable tilting module $T(\lambda)$ with
$[T(\lambda):L(\lambda)] =1$.

{\bf{Important convention:}
All weights 
$(a,b,c)$ will be denoted $(a-b,b-c)$.}

Normally we would label the highest weight modules by $\lambda \in
X^+$. However we don't want to have to
keep track of the degree of the representation. That is, we really
want to pretend we are looking at modules for $\SL_3(k)$, even though
such an object does not exist for the Dipper-Donkin quantisation, as
the determinant is not central. 
Since, however, we only need to consider polynomial modules and this
category splits up into a direct sum of homogeneous ones, we may
assume that we are always looking at modules of the same degree. 
Also we have the isomorphisms
$\nabla(a+d,b+d,c+d) \cong \nabla(a,b,c) \otimes D^{\otimes d}$, 
$L(a+d,b+d,c+d) \cong L(a,b,c) \otimes D^{\otimes d}$
and 
$T(a+d,b+d,c+d) \cong T(a,b,c) \otimes D^{\otimes d}$.
Thus we will label modules by the equivalent $\SL_3(k)$ weights. Thus
all the results in this paper will be in $\SL_3(k)$ notation (i.e. our
weights are in $\N^{\oplus 2}$). 
We may convert back by adding an appropriate power of the determinant
so that the modules all have the same degree.

\section{Preliminaries}\label{sect:prelim}

We first start off by noting the composition series of small induced
modules.
\begin{lem}\label{lem:small}
\begin{enumerate}
\item[(i)]{
Suppose $\lambda = (r,s)$ with $(r,s) \in \hat{C}$,
or $\lambda = (l-1,r)$ or $(r,l-1)$ with $0\le r \le l-1$.
Then $\nabla(\lambda)= L(\lambda)$. 
}
\item[(ii)]{
Suppose $\lambda = (l-s-2,l-r-2)$ with $(r,s) \in {C}$.
Then $\nabla(\lambda)$ has two composition factors 
with $L(\lambda)$ as its socle and
$L(r,s)$ as its head. 
}
\item[(iii)]{
Suppose $\lambda = l(1,0) +(r,s)$ with $(r,s) \in \hat{C}$.
Then $\nabla(\lambda)$ has two composition factors 
with $L(\lambda)$ as its socle and
$L(l-r-2,r+s+1)$ as its head. 
}
\item[(iv)]{
Suppose $\lambda = l(0,1) +(r,s)$ with $(r,s) \in \hat{C}$.
Then $\nabla(\lambda)$ has two composition factors 
with $L(\lambda)$ as its socle and
$L(r+s+1,l-r-2)$ as its head. 
}
\end{enumerate}
\end{lem}
This may be proved as in the classical case using Jantzen's sum
formula and translation functors.

%

%
\section{Translating the $\nabla_l$'s}
We start by considering the action of the translation functors
on the $\nabla_l$'s.

\begin{lem}
The translate of a $G$-module with a good $l$-filtration also
has a good $l$-filtration.
\end{lem}
\begin{proof}
This follows using the results of \cite{anderpfil} and the definition of
translation functors.
\end{proof}

We start by translating ``onto the walls''.
\begin{propn}\label{transontowall}
Let $\lambda$, $\mu \in \bar{C}$ such that $\mu$ belongs to the
closure of the facet containing $\lambda$. Let $w \in W_l$ with 
$w \cdot \lambda \in X^+$ and denote by $F$ the facet with $w\cdot
\lambda \in F$.
Then
$$
T_\lambda^\mu \nabla_l(w\cdot \lambda) \cong
\left\{\begin{array}{ll}
\nabla_l(w\cdot\mu), &\qquad\mbox{if $w\cdot\mu \in \hat{F}$,}\\
0, &\qquad\mbox{otherwise.}
\end{array}
\right.
$$
\end{propn}
\begin{proof}

Now by definition
$$
T_\lambda^\mu \nabla_l(w\cdot \lambda) \cong 
\pr_{\mu} (\nabla_l(w\cdot \lambda) \otimes L(\nu))$$
where $\nu$ is the unique element in $X^+ \cap W(\mu-\lambda)$,
(since $\nabla_l(w \cdot \lambda)$ is indecomposable).

We may use the tensor identity,  
\begin{align*}
T_\lambda^\mu \nabla_l(w\cdot \lambda) &\cong 
T_\lambda^\mu \Ind_{G_1B}^G \hat{L}(w\cdot \lambda)\\
&\cong 
\pr_{\mu} (\Ind_{G_1B}^G (\hat{L}(w\cdot \lambda)) \otimes L(\nu))\\
&\cong
\pr_{\mu} (\Ind_{G_1B}^G( \hat{L}(w\cdot \lambda) \otimes L(\nu)))\\
&\cong
\Ind_{G_1B}^G( \hat{\pr}_{\mu}( \hat{L}(w\cdot \lambda) \otimes L(\nu)))\\
&\cong
\Ind_{G_1B}^G( \hat{T}_{\lambda}^\mu( \hat{L}(w\cdot \lambda)))\\
&\cong
\left\{\begin{array}{ll}
\Ind_{G_1B}^G \hat{L}(w\cdot\mu), &\qquad\mbox{if $w\cdot\mu \in \hat{F}$,}\\
0, &\qquad\mbox{otherwise}
\end{array}
\right.\\
&\cong
\left\{\begin{array}{ll}
\nabla_l(w\cdot\mu), &\qquad\mbox{if $w\cdot\mu \in \hat{F}$,}\\
0, &\qquad\mbox{otherwise.}
\end{array}
\right.
\end{align*}
where we use $\hat{\ }$'s to distinguish modules and functors for
$G_1B$ from those for $G$. We also use the quantum version of \cite[II, remark 7.6
(1)]{jantz2} to identify $\hat{\pr_{\mu}}( - \otimes L(\nu))$ with the
translation functor $\hat{T}_{\lambda}^{\mu}$ on $\Mod(G_1B)$.
\end{proof}

\begin{rem}
We did not use the assumption that $G=q$-$\GL_3(k)$ or $\GL_3(k)$  thus
the above proposition is true for any quantum group  or linear
algebraic group $G$ where we
have the appropriate theory of $G_1B$-modules and translation functors.
\end{rem}

It will also be useful to know what happens when we translate back the
other way.
This is not as nice however and we will work it out on a case by case
basis.

\begin{propn}\label{propn:trans3}
Suppose $l \ge 3$.
Let $\lambda$, $\mu \in X^+$ with $\mu$ in the lower closure of the
alcove containing $\lambda$.
Then we have the following.
\begin{enumerate}
\item[(i)]
Suppose $\mu' = (l-1,r)$ with $0 \le r \le l-2$,
and $\lambda' = (a,b)$ with $0 \le a \le l-3$, and $0\le a+b \le l-3$.
Then
$T^{\lambda}_\mu \bnabla(\mu'')\frob \otimes L(\mu')$ 
has a good $l$-filtration with factors as shown.
$$
\xymatrix@R=15pt{
*=0{\bullet} \ar@{-}[d]^{\textstyle{
{\bnabla(\mu'')\frob \otimes L(l-a-2,a+b+1)}}} \\
*=0{\bullet} \ar@{-}[d]^{\textstyle{
{\bnabla(\mu''+(1,0))\frob \otimes L(\lambda')}}} \\
*=0{\bullet} \ar@{-}[d]^{\textstyle{
{\bnabla(\mu''+(-1,1))\frob \otimes L(\lambda')}}} \\
*=0{\bullet} \ar@{-}[d]^{\textstyle{
{\bnabla(\mu''+(0,-1))\frob \otimes L(\lambda')}}} \\
*=0{\bullet} \ar@{-}[d]^{\textstyle{
{\bnabla(\mu'')\frob \otimes L(l-a-b-3,a) }}}\\
*=0{\bullet} \ar@{-}[d]^{\textstyle{
{\bnabla(\mu'')\frob \otimes L(l-a-b-3,a) }}}\\
*=0{\bullet} \ar@{-}[d]^{\textstyle{
{\bnabla(\mu'')\frob \otimes L(l-a-2,a+b+1)}}}\\
*=0{\bullet} 
}
$$
\item[(ii)]{
Suppose $\mu' = (s,l-1)$ with $0 \le s \le l-2$,
and $\lambda' = (a,b)$ with $0 \le a \le l-3$, and $0\le a+b \le l-3$.
Then
$T^{\lambda}_\mu \bnabla(\mu'')\frob \otimes L(\mu')$ 
has a good $l$-filtration with factors as shown.
$$
\xymatrix@R=15pt{
*=0{\bullet} \ar@{-}[d]^{\textstyle{
{\bnabla(\mu'')\frob \otimes L(a+b+1,l-b-2)  }}}\\
*=0{\bullet} \ar@{-}[d]^{\textstyle{
{\bnabla(\mu''+(0,1))\frob \otimes L(\lambda') }}}\\
*=0{\bullet} \ar@{-}[d]^{\textstyle{
{\bnabla(\mu''+(1,-1))\frob \otimes L(\lambda')  }}}\\
*=0{\bullet} \ar@{-}[d]^{\textstyle{
{\bnabla(\mu''+(-1,0))\frob \otimes L(\lambda') }}}\\
*=0{\bullet} \ar@{-}[d]^{\textstyle{
{\bnabla(\mu'')\frob \otimes L(b,l-a-b-3) }}}\\
*=0{\bullet} \ar@{-}[d]^{\textstyle{
{\bnabla(\mu'')\frob \otimes L(a+b+1,l-b-2) }}}\\
*=0{\bullet} 
}
$$
}
\item[(iii)]{
Suppose $\mu' = (r,s)$ with $0 \le r \le l-2$ and $r+s=l-2$
and $\lambda'$ is in an up alcove.
Then
$T_\mu^{\lambda}\bnabla(\mu'')\frob \otimes L(\mu')$ 
has a good $l$-filtration with factors as shown.
$$
\xymatrix@R=15pt{
*=0{\bullet} \ar@{-}[d]^{\textstyle{
{\bnabla(\mu'')\frob \otimes L((l-2)\rho + w_0\lambda') }}}\\
*=0{\bullet} \ar@{-}[d]^{\textstyle{
{\bnabla(\mu'')\frob \otimes L(\lambda') }}}\\
*=0{\bullet} \ar@{-}[d]^{\textstyle{
{\bnabla(\mu'')\frob \otimes L((l-2)\rho + w_0\lambda') }}}\\
*=0{\bullet} 
}
$$
}
\end{enumerate}
\end{propn}
\begin{proof}
Case (i).
$$T_\mu^{\lambda}\bnabla(\mu'')\frob \otimes L(\mu')
\cong
\pr_{\lambda} \bnabla(\mu'')\frob \otimes \nabla(\mu')\otimes \nabla(\nu)
$$ 
We may use translation to assume that $\lambda'$ is such that $\nu =
(1,0)$.

Now $\nabla(l-1,r)\otimes \nabla(1,0)$ has a good filtration with
factors
(starting at the top)
$\nabla(l,r)$, $\nabla(l-2,r+1)$, $\nabla(l-1,r-1)$.

Thus the module
$\bnabla(\mu'')\frob \otimes \nabla(\mu')\otimes \nabla(\nu)$
has a filtration as shown,
$$
\xymatrix@R=15pt{
*=0{\bullet} \ar@{-}[d]^{\textstyle{
{\bnabla(\mu'')\frob \otimes L(l-2,r+1) }}}\\
*=0{\bullet} \ar@{-}[d]^{\textstyle{
{\bnabla(\mu'')\frob \otimes L(l,r)}}}\\
*=0{\bullet} \ar@{-}[d]^{\textstyle{
{\bnabla(\mu'')\frob \otimes L(l-r-3,0)  }}}\\
*=0{\bullet} \ar@{-}[d]^{\textstyle{
{\bnabla(\mu'')\frob \otimes L(l-2,r+1)}}}\\
*=0{\bullet} \ar@{-}[d]^{\textstyle{
{\bnabla(\mu'')\frob \otimes L(l-1,r-1)  }}}\\
*=0{\bullet} 
}
$$
using lemma \ref{lem:small}.
All the simples are $l$-restricted except for $L(l,r)$.

Now 
$$\bnabla(\mu'')\frob \otimes L(l,r)
\cong
 \bnabla(\mu'')\frob \otimes \bnabla(1,0)\frob \otimes L(0,r)$$
using Steinberg's tensor product theorem.
Also $\bnabla(\mu'') \otimes \bnabla(1,0)$ has a good filtration with
factors
(starting at the top)
$\bnabla(\mu''+(1,0))$, $\bnabla(\mu''+(-1,1))$,
$\bnabla(\mu''+(0,-1))$,
where the modules $\bnabla(\mu''+ (-1,1))$ and $\bnabla(\mu''
+(0,-1))$ are understood to be zero if the weight isn't dominant.

Now the weight $(l-1,r-1)$ is either not dominant or lies on a wall.
So after applying $\pr_{\lambda}$ to our filtration of 
$\bnabla(\mu'')\frob \otimes \nabla(\mu')\otimes \nabla(\nu)$
we get a module with good
$l$-filtration as shown.
$$
\xymatrix@R=15pt{
*=0{\bullet} \ar@{-}[d]^{\textstyle{
{\bnabla(\mu'')\frob \otimes L(l-2,r+1)}}}\\
*=0{\bullet} \ar@{-}[d]^{\textstyle{
{\bnabla(\mu''+(1,0))\frob \otimes L(0,r)}}}\\
*=0{\bullet} \ar@{-}[d]^{\textstyle{
{\bnabla(\mu''+(-1,1))\frob \otimes L(0,r)}}}\\
*=0{\bullet} \ar@{-}[d]^{\textstyle{
{\bnabla(\mu''+(0,-1))\frob \otimes L(0,r)}}}\\
*=0{\bullet} \ar@{-}[d]^{\textstyle{
{\bnabla(\mu'')\frob \otimes L(l-r-3,0)}}}\\
*=0{\bullet} \ar@{-}[d]^{\textstyle{
{\bnabla(\mu'')\frob \otimes L(l-2,r+1)}}}\\
*=0{\bullet} 
}
$$

We can use translation again to get the result as stated.

Case (ii).
This is the dual case to case (i).

Case (iii).
$$T_\mu^{\lambda}\bnabla(\mu'')\frob \otimes L(\mu')
\cong
\pr_{\lambda} \bnabla(\mu'')\frob \otimes \nabla(\mu')\otimes \nabla(\nu)
$$ 
We may use translation to assume that $\lambda'$ is such that $\nu =
(1,0)$.

Now $\nabla(r,s)\otimes \nabla(1,0)$ has a good filtration with
factors
(starting at the top)
$\nabla(r+1,s)$, $\nabla(r-1,s+1)$, $\nabla(r,s-1)$.

Thus the module
$\bnabla(\mu'')\frob \otimes \nabla(\mu')\otimes \nabla(\nu)$
has good $l$-filtration as shown,
$$
\xymatrix@R=15pt{
*=0{\bullet} \ar@{-}[d]^{\textstyle{
{\bnabla(\mu'')\frob \otimes L(r,s-1)}}}\\
*=0{\bullet} \ar@{-}[d]^{\textstyle{
{\bnabla(\mu'')\frob \otimes L(r+1,s)}}}\\
*=0{\bullet} \ar@{-}[d]^{\textstyle{
{\bnabla(\mu'')\frob \otimes L(r-1,s+1)}}}\\
*=0{\bullet} \ar@{-}[d]^{\textstyle{
{\bnabla(\mu'')\frob \otimes L(r,s-1)}}}\\
*=0{\bullet} 
}
$$
using lemma \ref{lem:small}.
The weight $(r-1,s+1)$ is either not dominant or lies on a wall, the
other simples are all $l$-restricted.
So after applying $\pr_{\lambda}$  we get a module with good
$l$-filtration
as above but without the $\bnabla(\mu'')\frob \otimes L(r-1,s+1)$.

We can use translation again to get the result as stated.
\end{proof}

A similar proof shows for $l=2$ that
\begin{propn}\label{propn:trans2}
Assume that $l=2$.
Let $\lambda$, $\mu \in X^+$ with $\mu$ in the lower closure of the
alcove for which $\lambda$ is in the upper closure.
Then we have the following.
\begin{enumerate}
\item[(i)]{
Suppose $\mu' = (1,0)$,
and $\lambda' = (0,0)$ 
Then
$T^{\lambda'}_\mu \bnabla(\mu'')\frob \otimes L(\mu')$ 
has a good $l$-filtration with factors as shown.
$$
\xymatrix@R=15pt{
*=0{\bullet} \ar@{-}[d]^{\textstyle{
{\bnabla(\mu'')\frob \otimes L(0,1)}}}\\
*=0{\bullet} \ar@{-}[d]^{\textstyle{
{\bnabla(\mu''+(1,0))\frob}}}\\
*=0{\bullet} \ar@{-}[d]^{\textstyle{
{\bnabla(\mu''+(-1,1))\frob}}}\\
*=0{\bullet} \ar@{-}[d]^{\textstyle{
{\bnabla(\mu''+(0,-1))\frob}}}\\
*=0{\bullet} \ar@{-}[d]^{\textstyle{
{\bnabla(\mu'')\frob \otimes L(0,1)}}}\\
*=0{\bullet} 
}
$$
}
\item[(ii)]{
Suppose $\mu' = (0,1)$,
and $\lambda' = (0,0)$.
Then
$T^{\lambda'}_\mu \bnabla(\mu'')\frob \otimes L(\mu')$ 
has a good $l$-filtration with factors as shown.
$$
\xymatrix@R=15pt{
*=0{\bullet} \ar@{-}[d]^{\textstyle{
{\bnabla(\mu'')\frob \otimes L(1,0)}}}\\
*=0{\bullet} \ar@{-}[d]^{\textstyle{
{\bnabla(\mu''+(0,1))\frob }}}\\
*=0{\bullet} \ar@{-}[d]^{\textstyle{
{\bnabla(\mu''+(1,-1))\frob}}}\\
*=0{\bullet} \ar@{-}[d]^{\textstyle{
{\bnabla(\mu''+(-1,0))\frob}}}\\
*=0{\bullet} \ar@{-}[d]^{\textstyle{
{\bnabla(\mu'')\frob \otimes L(1,0)}}}\\
*=0{\bullet} 
}
$$
}
\item[(iii)]{
Suppose $\mu' = (0,0)$ and $\lambda'=(1,0)$ or $(0,1)$.
Then
$$T_\mu^{\lambda'}\bnabla(\mu'')\frob \otimes L(\mu')
\cong
\bnabla(\mu'')\frob \otimes L(\lambda').
$$ 
}
\end{enumerate}
\end{propn}

We will also need.
\begin{propn}\label{propn:trans2b}
Assume that $l=2$.
Let $\lambda$, $\mu \in X^+$ with $\lambda$ and $\mu$ in the lower closure
of the same alcove but on different walls.
Then $\mu' = (1,0)$,
and $\lambda' = (0,1)$, 
or $\mu' = (0,1)$,
and $\lambda' = (1,0)$.
We have 
$$T^{\lambda'}_\mu \bnabla(\mu'')\frob \otimes L(\mu')
\cong \bnabla(\mu'')\frob.$$ 
\end{propn}
\begin{proof}
Now  $L(\mu') \otimes \nabla(1,0)$ has a good filtration with factors
$\nabla(1,1)$ and $\nabla(0)$. This splits as $\nabla(1,1)$ is the
Steinberg module. Thus 
\begin{align*}
\pr_{\lambda} \bnabla(\mu'')\frob \otimes \nabla(\mu')\otimes \nabla(1,0)
\cong 
& \, \bnabla(\mu'')\frob
\qedhere
\end{align*}
\end{proof}

\section{Characters}

Each $\nabla(\lambda)$ has an $l$-filtration (we may use the quantum
version of the argument of Jantzen \cite[3.13]{jandar}) but we would
like to know what the composition factors of $\hat{Z}(\lambda)$ are
for $\lambda \in X$. 

To do this we will work backwards -
and use the formula
\begin{equation}
 \ch \Ind_{G_1B}^G M = \sum_{\mu \in X} [M:\hat{L}(\mu)]\chi_l(\mu)
\label{chform}
\end{equation}
where $\chi_l(\mu) = \ch \nabla_l(\mu) = \chi(\mu'')\frob\phi(\mu')$ where
we put $\phi(\mu') = \ch L(\mu')$.
This is the quantum version of \cite[section 3]{jandar}. 

\begin{thm}\label{thm:lchar}
\begin{enumerate}
\item[(i)]{ Suppose 
$\lambda = l(a,b) + (l-1,l-1)$ with $(a,b) \in X^+$. Then
$\chi(\lambda)= \chi(a,b)\frob \phi(l-1,l-1)$.
}
\item[(ii)]{ Suppose
$\lambda = l(a,b) + (l-1,r)$ with $(a,b) \in X^+$ and $(l-1,r)\in
X_1$.
If we set $s=l-r-2$ then 
\begin{align*}
\chi(\lambda)&= \chi(a,b-1)\frob\phi(s,l-1) +
\chi(a+1,b-1)\frob\phi(r,s)\\
&\hspace{20pt}+ \chi(a-1,b)\frob \phi(r,s) + \chi(a,b)\frob\phi(l-1,r).
\end{align*}
These weights are depicted in Figure \ref{fig:walls}(a). 
\begin{figure}[ht]
\begin{center}
\epsfbox{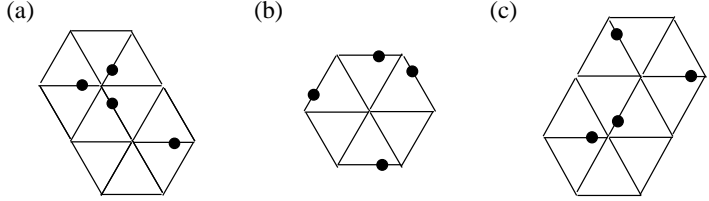}
\end{center}
\caption{\label{fig:walls}Diagram showing weights for $\lambda$ on (a)
  a right hand wall, (b) a left hand wall and (c) a horizontal wall}
\end{figure}
}
\item[(iii)]{ Suppose
$\lambda = l(a,b) + (s,l-1)$ with $(a,b) \in X^+$ and $(s,l-1)\in
X_1$.
If we set $r=l-s-2$ then 
\begin{align*}
\chi(\lambda)&= \chi(a-1,b)\frob\phi(l-1,r) +
\chi(a-1,b+1)\frob\phi(r,s)\\
&\hspace{20pt}+ \chi(a,b-1)\frob \phi(r,s) + \chi(a,b)\frob\phi(s,l-1).
\end{align*}
These weights are depicted in Figure \ref{fig:walls}(b). 
}
\item[(iv)]{ Suppose
$\lambda = l(a,b) + (r,s)$ with $(a,b) \in X^+$, $(r,s)\in X_1$ and
$r+s=l-2$. Then 
\begin{align*}
\chi(\lambda)&= 
\chi(a-1,b-1)\frob\phi(r,s) 
+ \chi(a,b-1)\frob\phi(l-1,r) \\
&\hspace{20pt} +\chi(a-1,b)\frob\phi(s,l-1) 
\chi(a,b)\frob\phi(r,s).
\end{align*}
These weights are depicted in Figure \ref{fig:walls}(c). 
}
\item[(v)]{ Suppose
$\lambda = l(a,b) + (r,s)$ with $(a,b) \in X^+$
and $(r,s)\in C$. We let 
$$\begin{array}{rlrl}\mu_1&=\lambda, 
&\mu_2&=(l a{+}r{+}s{+}1,l b{-}s{-}2),\\
\mu_3&=(la {+}l{-}r{-}s{-}3,l b{-}2 l{+}r), 
&\mu_4&=(l a{-}r{-}2,l b{+}r{+}s{+}1),\\
\mu_5&=(l a{-}2 l{+}s,l b {+}l{-}r{-}s{-}3), 
&\mu_6&=(l a{+}s,l b{-}r{-}s{-}3),\\
\mu_7&=(l a{-}l{+}r,l b{-}l{+}s), 
&\mu_8&=(la {-}r{-}s{-}3,l b{+}r), \\ 
\mu_9&=(l a{-}s{-}2,l b{-}r{-}2).& 
\end{array}$$
These weights are depicted in Figure \ref{fig:alc}(a), 
where the number corresponds to
the subscript of $\mu$.
\begin{figure}[ht]
\begin{center}
\epsfbox{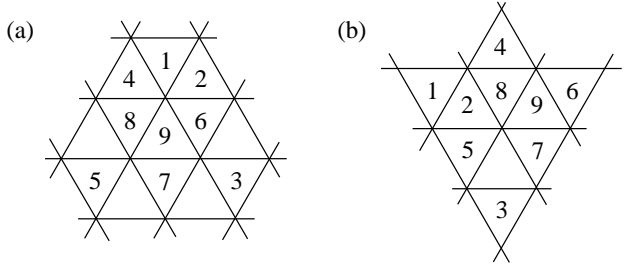}
\end{center}
\caption{\label{fig:alc}Diagram showing weights for $\lambda$ inside (a) a lower
alcove and (b) an upper alcove}
\end{figure}

Then $ \chi(\lambda) = \sum_{i=0}^{9} \chi_l(\mu_i)$.
}
\item[(vi)]{ Suppose
$\lambda = l(a,b) + (l-s-2,l-r-2)$ with $(a,b) \in X^+$, 
and $(r,s)\in C$. We let 
$$\begin{array}{rlrl}\mu_1&=(l a{-} l{+}s,l b {+}2
l{-}r{-}s{-}3),
 &\mu_2 &=(l a{-}r{-}2,l b{+}r{+}s{+}1),\\
\mu_3&=(l a{-}l{+}r,l b{-}l{+}s), &\mu_4&=\lambda,\\
\mu_5&=(la {-}r{-}s{-}3,l b{+}r), 
&\mu_6&=(la {+}2 l{-}r{-}s{-}3,l b{-}l{+}r),\\ 
\mu_7&=(l a{+}s,l b{-}r{-}s{-}3), &\mu_8&=(l a{+}r,l b{+}s) \\
\mu_9&=(l a{+}r{+}s{+}1,l b{-}s{-}2).
\end{array}$$
These weights are depicted in Figure \ref{fig:alc}(b). 

Then $ \chi(\lambda) = \sum_{i=0}^{9} \chi_l(\mu_i)$.
}
\end{enumerate}
\end{thm}
\begin{proof}
This is easily verified using induction and translation functors and
the previous propositions.

If $\lambda \in C$ then
$\chi_l(\mu_i) = 0$ for 
$2 \le i \le 9$.  For these $\mu_i$,  $\chi(\mu_i'') = 0$, as
$\mu_i''$ is fixed by one of the elements of $W$ under the dot action.
Thus 
$$ \sum_{i=0}^{9} \chi_l(\mu_i)
=  \chi_l(\mu_1) 
= \chi(\lambda)$$
using \ref{lem:small}.
We may use a similar argument for $\lambda \in \bar{C} \cap X^+$.

Now let $\lambda \in X^+$.
If $\lambda$ lies on a vertex then we have the well known result
that $\nabla(\lambda) \cong \bnabla(\lambda'')\frob \otimes
L(l-1,l-1)$ and thus have the required character formulae.

Suppose $\lambda$ lies on a wall and $l \ge 3$ - then we may translate
an induced module corresponding to a weight inside the alcove lying
below it ($\mu$ say) 
onto the wall. Since $T^\lambda_\mu \nabla(\mu) = \nabla(\lambda)$
we have 
$$\chi(\lambda) = \sum_{i} \ch( T^\lambda_\mu(\nabla_l(\mu_i)))$$
where $\mu_i$ are as in Figure \ref{fig:alc}.
We may now use proposition \ref{propn:trans3} to deduce the desired
character, noting that $\chi_l(\lambda_i)$ will be zero if one of the
parts of $\lambda_i''$ is $-1$.

If $\lambda$ lies inside an alcove (or lies on a wall and $l=2$) then
we may take a weight $\mu$ lying on a wall in the lower closure of (the
closure of) the alcove containing $\lambda$. 
Then  $\ch(T^\lambda_\mu\nabla(\mu)) = \ch(\nabla(\lambda)) +
\ch(\nabla(w\cdot \lambda))$, where $w$ is the unique reflection of
$W_l$ that fixes $\mu$. 
So
$$\chi(\lambda) = \sum_{i} \ch( T^\lambda_\mu(\nabla_l(\mu_i))) -
\chi(w\cdot \lambda)$$
where the $\mu_i$ will be (at most) four weights in the good
$l$-filtration of $\nabla(\mu)$.
The $\chi(w\cdot \lambda)$ is known by induction and the characters of
the translated $\nabla_l(\mu_i)$ may be deduced from proposition
\ref{propn:trans3} if $l \ge 3$ or propositions \ref{propn:trans2} and
\ref{propn:trans2b} if $l=2$.
Note that for generic $\mu$ and $l \ge 3$
the translate will have $6 + 6 + 2\times 3
= 18$ factors as one would expect from adding the factors of
$\nabla(\lambda)$ and $\nabla(w \cdot \lambda)$.
For generic $\mu$ and $l=2$ then the translate has
$5+1+1+1=8$ factors.

Also note that if $\lambda$ is in a down alcove and is right on the edge of the
dominant region (that is $\lambda'' =(a,0)$ or $(0,a)$ for some $a \in
\N$), then $\chi_l(\mu_3) = - \chi_l(\mu_8)$ so these cancel in the
sum.
\end{proof}

\begin{cor}
We have the following characters for $\hat{Z}(l\lambda'' + \lambda')$
with $\lambda''
\in X$ and $\lambda' \in X_1$.
\begin{enumerate}
\item[(i)]{Suppose $\lambda' = (l-1,l-1)$, then
$$\ch \hat{Z}(l\lambda'' + (l-1,l-1))
= \ch \hat{L}(l\lambda''+ (l-1,l-1)).$$
}
\item[(ii)]{ Suppose $\lambda'=(l-1,r)$ with $0 \le r \le l-2$, then
\begin{align*}
\ch \hat{Z}(l\lambda'' + (l-1,r)) 
&= \ch \hat{L}((l-1,r)+l \lambda'' )+  
\ch \hat{L}((r-l,s)+ l \lambda'') \\
&\hspace{20pt}+\ch \hat{L}((r+l,s-l)+ l \lambda'')
+\ch \hat{L}((s,-1)+ l \lambda'').
\end{align*}
}
\item[(iii)]{ Suppose $\lambda'=(s,l-1)$ with $0 \le s \le l-2$, then
\begin{align*}
\ch \hat{Z}(l\lambda'' + (s,l-1)) 
&= \ch \hat{L}((s,l-1)+l \lambda'' )+  
\ch \hat{L}((r,s-l)+ l \lambda'') \\
&\hspace{20pt}+\ch \hat{L}((r-l,s+l)+ l \lambda'')
+\ch \hat{L}((-1,r)+ l \lambda'').
\end{align*}
}
\item[(iv)]{ Suppose $\lambda'=(r,s)$ with  $0 \le r \le l-2$ and $r+s=l-2$,
then
\begin{align*}
\ch \hat{Z}(l\lambda'' + (r,s)) 
&= \ch \hat{L}((r,s)+l \lambda'' )+  
\ch \hat{L}((s-l,l-1)+ l \lambda'') \\
&\hspace{20pt}+\ch \hat{L}((l-1,r-l)+ l \lambda'')
+\ch \hat{L}((r-l,s-l)+ l \lambda'').
\end{align*}
}
\item[(v)]{ Suppose $\lambda' = (r,s) \in C$, then
$$\ch \hat{Z}(l\lambda'' + (r,s)) 
= \sum_{i} \ch \hat{L}(\mu_i)$$ 
where the $\mu_i$ are as in Figure~\ref{fig:alc}(a).
}
\item[(vi)]{Suppose $\lambda'= (l-s-2,l-r-2)$ with $(r,s) \in C$, then
$$\ch \hat{Z}(l\lambda'' + (l-s-2,l-r-2)) 
= \sum_{i} \ch \hat{L}(\mu_i)$$
where the $\mu_i$ are as in Figure~\ref{fig:alc}(b).
}
\end{enumerate}
\end{cor}
\begin{proof}
We have $ \Ind_{G_1B}^G \hat{Z}(\lambda) \cong \nabla(\lambda)$ so
this follows using the character formula \eqref{chform}, 
the previous theorem and the identity 
$$ \hat{Z}(\lambda' + l \lambda'') \cong 
 \hat{Z}(\lambda') \otimes k_{l \lambda''} $$
which is the quantum version of \cite[II 9.2 (5)]{jantz2}. The quantum
 result follows as in the classical case.
\end{proof}

\section{Extensions for simple modules}
We will need to be able to work out the $G_1B$ extensions between
simple modules for $G_1B$. To do this we will need to generalise the
extension results of Yehia \cite{yehia}. We will essentially reproduce
his proofs but in the quantum case, as the reference is not widely
accessible.

\begin{lem}
Let $\lambda \in X_1$ then $L(\lambda) \otimes \St$ has a good
filtration.
\end{lem}
\begin{proof}
If $l \ge4=2h-2$ then this is the quantum
version of \cite[2.5 corollary]{anderpfil}.  

If $\lambda$ is not in an up alcove then $L(\lambda) \cong
\nabla(\lambda)$ and we are done by \cite{mathieu} and \cite[corollary
5.14]{andpolwen}.

So the only case left is if $l=3$ and $\lambda = (1,1)$.
But now $\ch(L(1,1))= \ch(\nabla(1,1)) - \ch(\nabla(0,0)) =
e(1,1) + e(2,-1) +e(1,-2) + e(-1,-1) + e(-2,1) + e(-1,2) + e(0,0)$.
So all the weights of $L(1,1)\vert_B \otimes k_{(2,2)}$ are dominant and
so $\Ind_B^G L(1,1)\otimes k_{(2,2)} = L(1,1)\otimes \St$ has a good
filtration.
\end{proof}

\begin{propn}
Let $\lambda \in X_1$. 
There is an indecomposable $G$-module $Q(\lambda)$ which restricts to the $G_1$
injective hull of $L(\lambda)$ and this module is a tilting module for
$G$.
Moreover $Q(\lambda)$ is the tilting module $T(2(l-1)\rho
+w_0\lambda)$ and this module is a direct summand of the module 
$L((l-1)\rho +w_0\lambda) \otimes \St$.
\end{propn}
\begin{proof}
If $l \ge 4$ then this is the result \cite[proposition 5.7]{anderpfil}.

Let $\nu = (l-1)\rho + w_0\lambda \in X_1$.
So $L(\nu)^* \cong L((l-1)\rho -\lambda)$.
If $l \le 3$ and $\lambda$ lies on a left or right hand wall
then the tilting module
$T(2(l-1)\rho +w_0\lambda)$ 
is $T_{(1,1)}^\lambda \St \cong \pr_{\lambda} L(\nu)\otimes \St$. This then has
simple $G$-socle $L(\lambda)$ and is injective as a $G_1$-module. 
Let $\mu \in X_1$.
We have
$$
\Hom_{G_1}(L(\mu), L(\nu)\otimes \St)
\cong
\Hom_{G_1}(L(\mu)\otimes L((l-1)\rho -\lambda), \St)
$$
and the latter group has dimension 
$[L(\mu)\otimes L((l-1)\rho -\lambda): \St]_{G_1}$, the $G_1$
composition multiplicity of $\St$ in $L(\mu)\otimes L((l-1)\rho
-\lambda)$, 
as $\St$ is the $G_1$ injective hull of $\St$.
We may check that 
$$[L(\mu)\otimes L((l-1)\rho -\lambda): \St]_{G_1}
\cong
\left\{\begin{array}{ll}
1, &\qquad\mbox{if $\mu = \lambda$,}\\
0, &\qquad\mbox{otherwise.}
\end{array}
\right.
$$ 
Thus $L(\nu)\otimes \St \cong T(2(l-1)\rho +w_0\lambda$ and is the
$G_1$ injective hull of $L(\lambda)$.

If $\lambda =(0,0)$ and $l= 2$ then 
$\nu = (1,1)$.
We may check that 
$$
[L(\mu)\otimes \St: \St]_{G_1}
\cong
\left\{\begin{array}{ll}
1, &\qquad\mbox{if $\mu = (0,0)$,}\\
3, &\qquad\mbox{if $\mu = (1,1)$,}\\
0, &\qquad\mbox{otherwise.}
\end{array}
\right.
$$
Thus the module $\St \otimes \St$ is the direct sum of three copies of
the Steinberg module and one copy of the $G_1$ injective hull of
$L(0,0)$ which is $\pr_{(0,0)}(\St \otimes \St) \cong T(2,2)$.

If $\lambda=(1,1)$ and $l=3$ then the
translate $T_{(2,2)}^{\lambda} \St = \pr_{\lambda} L(1,1)\otimes \St$.
We may check that 
$$
[L(\mu)\otimes L(1,1): \St]_{G_1}
\cong
\left\{\begin{array}{ll}
1, &\qquad\mbox{if $\mu = (1,1)$ or $\mu=(2,2)$,}\\
0, &\qquad\mbox{otherwise.}
\end{array}
\right.
$$
Thus $L(1,1)\otimes \St$ is the direct sum of the Steinberg module and
the $G_1$-injective hull of $L(1,1)$ which is
 $\pr_{(1,1)}(L(1,1) \otimes \St) \cong T(3,3)$.

We may now get the $G_1$ injective hull of $L(0,1)$ or $L(1,0)$ by 
translating the
$T(3,3)$ onto the wall. This translate is $T(3,4)$ or $T(4,3)$
respectively. 
A similar argument to above shows that this module is injective as a
$G_1$ module and has $G_1$ socle $L(0,2)$ or $L(2,0)$ respectively.
Also the module $L(l-2,l-1)\otimes \St$ is a tilting module, a
character calculation shows that $T(3,4)$ is a direct summand of this
module.

If $\lambda = (0,0)$ and $l=3$  
then the translate  $ T_{(0,1)}^{(0,0)} T(4,3)=\pr_{(0,0)}L(0,1) \otimes
T(4,3)$ is injective as a
$G_1$-mod as it is a direct summand of a tensor product of an
injective $G_1$-module.
As a $G$-module $L(0,1) \otimes T(4,3)$ is isomorphic to $T(5,5) \oplus
T(5,2)\oplus T(5,2)$. 
We have
$$
\Hom_{G_1}(L(\mu), L(0,1)\otimes T(4,3))
\cong
\Hom_{G_1}(L(\mu)\otimes L(1,0), T(4,3))
$$
the latter group has dimension equal to the $G_1$ composition
multiplicity of $L(1,0)$ in $L(\mu) \otimes L(1,0)$ as $T(4,3)$ is the
$G_1$ injective hull of $L(1,0)$.
We may check that for $\mu \in X_1$ 
$$
[L(\mu)\otimes L(1,0): L(1,0)]_{G_1}
\cong
\left\{\begin{array}{ll}
1, &\qquad\mbox{if $\mu = (0,0)$}\\
6, &\qquad\mbox{if $\mu = (2,2)$}\\
0, &\qquad\mbox{otherwise.}
\end{array}
\right.
$$
Since $T(5,2)\cong \nabla(1,0)\frob \otimes \St$ we have
$\Hom_{G_1}(L(\mu), T(5,2))
\cong \nabla(1,0)\frob$ if $\mu = (2,2)$ and zero otherwise.
Thus $\Hom_{G_1}(L(\mu), T(4,3))$ is $k$ if $\mu = (0,0)$ and zero
otherwise and hence $T(4,3)$ is the $G_1$-injective hull of $L(0,0)$.

For $l=3$ the module $\St \otimes \St$ is a tilting module and it has
summands $T(4,4)$, $T(3,3)$, $T(5,2)$, $T(2,5)$ and three copies of
the Steinberg module, by characters.
\end{proof}

\begin{cor}\label{smallhead}
The $G$-head of $\nabla(2(l-1)\rho +w_0\lambda)$ is simple and is
isomorphic to $L(\lambda)$.
\end{cor}
\begin{proof}
We have that $\hd(\nabla(2(l-1)\rho +w_0\lambda) \subseteq \hd
T(2(l-1)\rho +w_0\lambda) \cong L(\lambda)$ by the previous  
proposition.
\end{proof}

The following four results follow as in the classical case
\cite[4.8-4.11]{parker1}, see also \cite{jandar}. 
\begin{cor}
If $\lambda \in X^+$ and $\mu \in X_1$ then $\bnabla(\lambda)\frob
\otimes T(2(l-1)\rho +w_0\mu)$ has a good filtration.
\end{cor}

\begin{cor}
If $\lambda \in X^+$ and $\mu \in X_1$ then $\nabla(l\lambda
+2(l-1)\rho +w_0\mu)$ is a quotient
of $\bnabla(\lambda)\frob
\otimes T(2(l-1)\rho +w_0\mu)$ and
$\nabla(l\lambda+\mu)$ as a submodule. 
\end{cor}

\begin{cor}
For all $\lambda \in X^+$ and $\mu \in X_1$ we have
$$\hd_{G_1} \nabla(l\lambda +2(l-1)\rho +w_0\mu)
\cong \bnabla(\lambda)\frob \otimes L(\mu) $$
and
$$\soc_{G_1} \nabla(l\lambda+ \mu)
\cong \bnabla(\lambda)\frob \otimes L(\mu). $$
\end{cor}

\begin{cor}
For all $\lambda \in X^+$ the module $\nabla(\lambda)$ has simple
head.
\end{cor}

To determine $\Ext_{G_1}^1(L(\mu), L(\lambda))$ we need to determine
$\Ext^1_G(L(\mu), L(\lambda))$ for small $\mu$ and $\lambda$.
``Small'' in this case means that $\lambda \le 2(l-1)\rho$ and $\mu
\in X_1$.

The idea is to use the quantum version of the short exact sequence
\cite{donkext}
\begin{multline} \label{five}
0 \to
\Ext^1_{G/G_1}(k, \Hom_{G_1}(L(\mu), L(\lambda)))
\to \Ext_{G}^1(L(\mu), L(\lambda))\\
\to
\Hom_{G/G_1}(k, \Ext^1_{G_1}(L(\mu), L(\lambda)))
\to 0
\end{multline}

Note that
$$\Hom_{G/G_1}(k, \Ext^1_{G_1}(L(\mu), L(\lambda)))
\cong\Hom_{G/G_1}(L(\mu'')\frob, \Ext^1_{G_1}(L(\mu'), L(\lambda))).
$$
Also $\Ext^1_{G_1}(L(\mu'), L(\lambda))\cong 
\Hom_{G_1}(L(\mu'), Q(\lambda)/L(\lambda))$ 
so determining $\Ext^1_G(L(\mu), L(\lambda))$ for enough $\mu$
determines the $G_1$
socle of  $Q(\lambda)/L(\lambda)$ which in turn determines
$\Ext^1_{G_1}(L(\mu'), L(\lambda))$.
We thus only need to calculate the $\Ext$ groups for $\mu$ a
composition factor of $Q(\lambda)$. I.e., it is enough to determine the
$\Ext$'s for $\mu \le 2(l-1)\rho$ and $\mu$ in the same block as
$\lambda$.

\begin{lem}\label{smallextwall}
Suppose $0 \le r \le l-2$ and $r+s=l-2$ then
$$\Ext^1_G(L(r,s),L(2l-1,r)) \cong 
\Ext^1_G(L(r,s),L(s,2l-1))\cong k.
$$
If $l \ne 3$ then 
$$\Ext^1_G(L(r,s), L(l+r,l+s)) \cong 0.
$$
If $l=3$ then 
$$\Ext^1_G(L(r,s), L(l+r,l+s)) \cong k.$$
\end{lem}
\begin{proof}
Since if $\mu \not >  \lambda$ we have
$\Ext^1_G(L(\mu), L(\lambda)) \cong \Hom_{G}(L(\mu),
\nabla(\lambda)/L(\lambda))$,
this lemma will follow if we know what the socle of
$\nabla(\lambda)/L(\lambda)$ is.

Now if $\lambda = (2l-1,r)$ or $(s,2l-1)$ then $\nabla(\lambda)$ only has two
composition factors $L(\lambda)$ and $L(r,s)$. Thus
$\nabla(\lambda)/L(\lambda) \cong L(r,s)$ and the result follows.

If $\lambda =(l+r, l+s)$ and $l\ne 3$ then $\nabla(\lambda)$ has four
composition factors: $L(\lambda)$, $L(l-1,r)$, $L(l-1,s)$ and
$L(r,s)$.
The previous corollary says that $L(r,s)$ is the head of
$\nabla(\lambda)$.
We also know that $\Ext^1_G(L(s,l-1), L(l-1,r)) \cong 
\Ext^1_G(L(l-1,r), L(s,l-1)) \cong 0 $ thus the socle of
$\nabla(\lambda)/L(\lambda)$ is $L(s,l-1) \oplus L(l-1,r)$.
Thus 
$\Ext^1_G(L(r,s), L(l+r,l+s) \cong 0$. 

If $l=3$ then
$[\nabla(l+r,l+s):L(r,s)]=2$. The module $\nabla(l+r,l+s)$ has simple
head $L(r,s)$. Since $\nabla(l+r,l+s)$ has five composition factors in
total and is indecomposable the multiplicity of $L(r,s)$ in socle of
$\nabla(l+r,l+s)/L(l+r,l+s)$ is at most one. Thus the dimension of 
$\Ext^1_G(L(r,s), L(l+r,l+s))$ is at most one. But there is at least
one non-split extension - it is the indecomposable module 
$\bnabla(1,1)\frob \otimes L(r,s)$.
\end{proof}

We similarly get:
\begin{lem}\label{smallextwall2}
Suppose $0 \le r \le l-2$ and $r+s=l-2$ then
$$\Ext^1_G(L(l-1,r),L(r,l+s)) \cong 
\Ext^1_G(L(s,l-1),L(l+r,s)\cong k.
$$
and
$$\Ext^1_G(L(l-1,r), L(l+s,l-1)) \cong 
\Ext^1_G(L(s,l-1), L(l-1,l+r)) \cong 
0.$$
\end{lem}

\begin{lem}\label{smallextupalc}
Suppose $0 \le r\le l-3$ and $0 \le r+s \le l-3$ then
$$\Ext^1_G(L(l-s-2,l-r-2),L(\nu))\cong k$$
if $\nu \in \{(r,s), (l+s,l-r-s-3), (l-r-s-3,l+r)\}$
and
$$\Ext^1_G(L(l-s-2,l-r-2),L(\nu))\cong 0$$
if $\nu \in \{(l-s-2, l-r-2), (2l-s-2,l-r-2),
(l-s-2,2l-r-2),(l+r,l+s)\}$.
\end{lem}
\begin{proof}
The result for the first $\Ext$ group follows from the fact 
that there are only two composition factors of
$\nabla(\nu)$ and $\nabla(l-s-2,l-r-2)$.

For the second $\Ext$ group 
we use that fact that $\nabla(\nu)$ (if $\nu \ne (l-s-2,l-r-2)$) has
simple head $L(l-s-2,l-r-2)$ and this is the only occurrence of this
simple module in $\nabla(\nu)$.
We may deduce that $\nabla(\nu)$ has simple head $L(l-s-2,l-r-2)$ by
either using corollary \ref{smallhead} or by translating an induced
module off the wall.
\end{proof}

\begin{lem}\label{smallextdownalc}
Suppose $0 \le r\le l-3$ and $0 \le r+s \le l-3$ then
$$\Ext^1_G(L(r,s),L(\nu))\cong k$$
if $\nu \in \{(l-s-2, l-r-2), (l-r-2,l+r+s+1),
(l+r+s+1,l-s-2)\}$
and
$$\Ext^1_G(L(r,s),L(\nu))\cong 0$$
if $\nu \in \{(r,s), (l+s,l-r-s-3), (l-r-s-3,l+r),
(s,3l-r-s-3), (3l-r-s-3,r), (2l-s-2, 2l-r-2) \}$.
If $l \ne 3$ then 
$$\Ext^1_G(L(r,s), L(l+r,l+s)) \cong 0.
$$
If $l=3$ then 
$$\Ext^1_G(L(r,s), L(l+r,l+s)) \cong k.$$
\end{lem}
\begin{proof}
We first observe that $\nabla(r+s+1,2l-s-2)$ is a quotient of
$\nabla(l+r,l+s)$ (and dually so is $\nabla(2l-r-2, r+s+1)$).
These modules all have the same simple head --- namely $L(l-s-2,l-r-2)$.
Also there is a unique homomorphism from $\nabla(l+r,l+s)$ to
$\nabla(r+s+1,2l-s-2)$. (Quantum version \cite[section 7]{andpolwen} of \cite[II,
  7.19(d)]{jantz2}.)
Since this homomorphism must be non-zero on the head of $\nabla(l+r,l+s)$ and this
head is the same as the head of $\nabla(r+s+1, 2l-s-2)$ and this
simple module only occurs once in $\nabla(r+s+1, 2l-s-2)$ this map
must be onto.

Thus by considering the composition factors of the kernel of this
homomorphism,
the socle of the quotient $\nabla(l+r,l+s)/L(l+r,l+s)$ is contained in 
$L(2l-r-2,r+s+1) \oplus L(r+s+1,2l-r-2)$ if $l\ne 3$ and 
$L(4,1) \oplus L(1,4) \oplus L(0,0)$ if $l= 3$.

Thus $\Ext_G^1(L(r,s), L(l+r,l+s)$ is zero if $l\ne 3$.
If $l=3$ then $\Ext_G^1(L(0,0),L(3,3))$ is at most one-dimensional. 
But there is a non-split extension --- namely the module
$\nabla(1,1)\frob$.

If $\nu \in \{(r,s), (l+s,l-r-s-3), (l-r-s-3,l+r)\}$
then $L(r,s)$ is not a composition factor of $\nabla(\nu)$ so
$\Ext_G^1(L(r,s), L(\nu)) \cong 0$.

We may now deduce that the socle of the quotient 
$\nabla(r+s+1, 2l-s-2)/L(r+s+1,2l-s-2)$ is 
$L(r,s)\oplus L(l+s,l-r-s-3) \oplus L(l-r-s-3,l+r)$ as these cannot
extend each other and the only other composition factor of
$\nabla(r+s+1,2l-s-2)$ is its head $L(l-s-2,l-r-2)$.
Thus $\Ext_G^1(L(r,s), L(r+s+1, 2l-s-2)) \cong k$.
Dually we have
$\Ext_G^1(L(r,s), L(2l-r-2, r+s+1)) \cong k$.

If $\nu = (l-s-2,l-r-2)$ then this extension is the module
$\nabla(l-s-2,l-r-2)$.

If $\nu \in (s,3l-r-s-3), (3l-r-s-3,r), (2l-s-2, 2l-r-2) \}$ 
and $l\ne 3$ then $L(r,s)$ is the head of $\nabla(\nu)$. Since
$\nabla(\nu)$ has both simple head and socle and has at least three
composition factors and $L(r,s)$ occurs with multiplicity one, it
cannot be in the socle of the quotient $\nabla(\nu)/L(\nu)$ thus 
$\Ext_G^1(L(r,s), L(\nu))$ is zero.

If $l=3$ the only case that the above paragraph does not work is for
$\nu = (4,4)$ when $L(0,0)$ occurs with multiplicity two.
If $\Ext_G^1(L(0,0), \bar{L}(1,1)\frob \otimes L(1,1))$ is non-zero then 
using the five term exact sequence $\bar{L}(1,1)\frob$ must be a composition
factor of $\Ext^1_{G_1}(L(0,0), L(1,1))$. The following lemma will show
that this is not the case and so $\Ext^1_G(L(0,0), L(4,4))$ is zero.
\end{proof}

\begin{lem}
If $l=3$ then 
$$\Ext_{G_1}^1(L(0,0), L(1,1)) \cong \bnabla(1,0)\frob \oplus
\bnabla(0,1)\frob \oplus k.$$
\end{lem}
\begin{proof}
The $G_1$ injective hull of $L(1,1)$ is $T(3,3)$.
We apply $\Hom_{G_1}(k, -)$ to the \ses
$$ 0 \to L(1,1) \to T(3,3) \to Q \to 0$$
to get
$$
0 \to \Hom_{G_1}(k, L(1,1))
\to \Hom_{G_1}(k, T(3,3))
\to \Hom_{G_1}(k, Q)
\to \Ext^1_{G_1}(k, L(1,1))
\to 0
$$
The first two $\Hom$ groups are zero so the last two groups are
isomorphic.
Thus $Q^{G_1} \cong  \Ext^1_{G_1}(k,\break L(1,1))$.

Now the $G_1$ fixed points of $Q$ are contained in the $G_1$ fixed
points of the
induced modules appearing in a good filtration of $T(3,3)/\nabla(1,1)$ together
with the $G_1$ fixed points of $\nabla(1,1)/L(1,1)$.
We thus have
$$
Q^{G_1}  \subseteq
k \oplus \bnabla(1,0)\frob \oplus \bnabla(0,1)\frob 
\oplus \bar{L}(1,1) \frob
$$
But $L(1,1)\frob $ can't be in the $G_1$ socle of $Q$ as then it would
also be in the $G_1$ head of the $Q^*$. 
The $G_1$ head of $Q^*$ is contained in the $G_1$ heads of the 
induced modules appearing in a good filtration of $T(3,3)$ as $T(3,3)$
is self dual.
Thus 
$$\hd_{G_1}(Q^*) \subseteq L(1,1) ^{\oplus 5} \oplus k^{\oplus
  2}\oplus
\bnabla(1,0)\frob \oplus \bnabla(0,1)\frob.$$

Hence 
$$
Q^{G_1}  \subseteq
k \oplus \bnabla(1,0)\frob \oplus \bnabla(0,1)\frob.
$$
We now observe from the good filtration of $T(3,3)$ that 
$\bnabla(1,0)\frob \oplus \bnabla(0,1)\frob$ must occur directly above
$k$ in $T(3,3)/L(1,1)$. The previous lemma tells us that $k$ cannot
extend
either $\bnabla(1,0)\frob$ nor $\bnabla(0,1)\frob$ so this is indeed the
$G_1$ fixed points of $Q$.
\end{proof}

We may now prove the following.

\begin{thm} \label{thm:yehia}
The $\Ext^1_{G_1}\bigl(L(\alpha),L(\beta)\bigr)$ for $\alpha$, $\beta \in
X_1$ are given by the following tables.
(i)\ 
For $(r,s) \in X_1$ with $r+s=l-2$, we have\\
\medskip

\noindent\begin{tabular}{c|ccc} 
$\alpha \downarrow$, $\beta \rightarrow$&$(r,s)$& $(l-1,r)$ 
& $(s,l-1)$ \\\hline 
$(r,s)$   &$0$    &$\bnabla(0,1)\frob$&$\bnabla(1,0)\frob$ \\
$(l-1,r)$ &$\bnabla(1,0)\frob$& $0$     & $0$ \\
$(s,l-1)$ &$\bnabla(0,1)\frob$& $0$     & $0$ \\
\end{tabular}\\
\medskip

\noindent(ii)\ 
For $(r,s) \in C$ and $ l \ge 4$, the only non-zero entries we have\\
\medskip

\noindent\begin{tabular}{c|ccc} 
$\alpha \downarrow$, $\beta \rightarrow$ &$(l-s-2,l-r-2)$   
& $(r+s+1,l-s-2)$ & $(l-r-2,r+s+1)$ \\ \hline
$(r,s)$   &$k$    &$\bnabla(0,1)\frob$&$\bnabla(1,0)\frob$ \\
\end{tabular}\\

\medskip

\noindent\begin{tabular}{c|ccc} 
$\alpha \downarrow$, $\beta \rightarrow$&$(r,s)$   
& $(s,l-r-s-3)$ & $(l-r-s-3,r)$ \\ \hline
$(l-s-2,l-r-2)$   &$k$               
&$\bnabla(0,1)\frob$&$\bnabla(1,0)\frob$ \\
\end{tabular}\\
\medskip

\noindent If $l=3$ then all the entries in the two tables above are replaced by
$k \oplus \bnabla(0,1)\frob\oplus \bnabla(1,0)\frob$.
\end{thm}
\begin{proof}
We use the sequence \eqref{five} and the previous results to show that
the $\Ext^1_{G_1}$ are as described. We have to argue as in the
previous lemma to do the case $l=3$.
\end{proof}

To now determine $\Ext^1_G(L(\mu), L(\lambda))$ for $\mu$ and $\lambda
\in X^+$ we need to know the $G_1$ socle of the tensor products
$L(1,0)\otimes L(\lambda)$ and $L(0,1)\otimes L(\lambda)$ for $\lambda \in X_1$. 
We essentially determined the tensor product in the proofs of propositions 
\ref{transontowall}, \ref{propn:trans3}, \ref{propn:trans2} and  
\ref{propn:trans2b}. We just need to determine the socles of these
tensor products.
These are not hard to compute using translation functors and follow
exactly as in the classical case so we will
just state the result.

\begin{propn}
The $G_1$ socle of the tensor product $L(1,0)\otimes L(\lambda)$ for
$\lambda \in X_1$ is the same as its $G$ socle and is given by the
following table. 
\\
\medskip

\noindent
\begin{tabular}{c|c|c} 
$l$ &$\lambda$ & $\soc_G L(1,0) \otimes L(\lambda)$\\ \hline
all $l$&$(0,0)$  & $L(1,0)$\\
$l\ge4$&$(0,s)$, $1\le s \le l-3$ & $L(1,s) \oplus L(0,s-1)$\\
$l\ge3$&$(0,l-2)$ & $L(0,l-3)$\\
$l\ge4$&$(r,s)$, $1\le r \le l-3$ and $r+s=l-2$ & $L(r,s-1) \oplus L(r-1,s+1)$\\
$l\ge3$&$(r,0)$, $1\le r \le l-2$ & $L(r+1,0) \oplus L(r-1,1)$\\
$l\ge4$&$(r,s)$ deep inside $C$ 
             & $L(r+1,s) \oplus L(r-1,s+1) \oplus L(r,s-1)$\\
all $l$&$(0,l-1)$ & $L(0,1,l-1) \oplus L(0,l-2)$\\
$l\ge3$&$(r,l-1)$, $1\le r \le l-2$ & $L(r+1,l-1) \oplus L(r,l-2)$\\
all $l$&$(l-1,l-1)$ & $L(l-1,l-2)$\\
$l\ge3$&$(1,l-2)$ & $L(2,l-2)\oplus L(0,l-1)$\\
$l\ge4$&$(r,l-2)$, $2\le r \le l-2$ & $L(r+1,l-2) \oplus L(r,l-3) \oplus
  L(r-1,l-3)$\\
$l\ge4$&$(r,s)$, $2\le r \le l-3$ and  $r+s=l-1$ & $L(r+1,s) \oplus
  L(r-1,s+1)$\\
$l\ge4$&$(l-2,1)$ & $L(l-1,1)\oplus L(l-3,2)$\\
all $l$&$(l-1,0)$  & $L(l-2,1)$\\
$l\ge3$&$(l-1,s)$, $1\le s \le l-2$ & $L(l-2,s+1) \oplus L(l-1,s-1)$\\
$l\ge4$&$(l-2,s)$, $2\le s \le l-2$ & $L(l-1,s) \oplus L(l-2,s-1) \oplus
  L(l-3,s+1)$\\
$l\ge4$&$(r,s)$ deep inside upper alcove 
             & $L(r+1,s) \oplus L(r-1,s+1) \oplus L(r,s-1)$\\
\end{tabular}
\end{propn}

We may use the dual of the above table to determine $L(0,1) \otimes
L(\lambda)$ for $\lambda \in X_1$.

\begin{cor}
Let $\lambda \in X^+$. 
Then
$$\soc_G L(1,0) \otimes L(\lambda) = (\soc_G L(1,0) \otimes
L(\lambda'))\otimes L(\lambda '')\frob$$
and 
$$\soc_G L(0,1) \otimes L(\lambda) = (\soc_G L(0,1) \otimes
L(\lambda'))\otimes L(\lambda '')\frob$$
\end{cor}
\begin{proof}
We have
$$\soc_G L(1,0) \otimes L(\lambda) = \soc_G(\soc_{G_1}(L(1,0) \otimes
L(\lambda'))\otimes L(\lambda '')\frob),$$
but the $G_1$ socle of $L(1,0) \otimes L(\lambda')$ 
is the same as its $G$ socle.
Steinberg's tensor product theorem then tells us that that
$\soc_{G_1}(L(1,0) \otimes L(\lambda'))\otimes L(\lambda'')\frob$ is
semi-simple as a $G$-module.
\end{proof}

We may now deduce the following theorem.
\begin{thm}
Let $\mu$, $\lambda \in X^+$. 
If $\mu' = \lambda'$ then
$\Ext_G^1(L(\mu), L(\lambda)) \cong \Ext_G^1(L(\mu''), L(\lambda'')$.

If $\mu' \ne \lambda'$ then
$\Ext_G^1(L(\mu), L(\lambda)) \cong
\Hom_{G}(L(\mu''),\Ext^1_{G_1}(L(\mu'),
L(\lambda'))^{(-1)}\otimes  L(\lambda'')$.

We have
$\dim \Ext_G^1(L(\mu), L(\lambda) \le 1$.
\end{thm}
\begin{proof}
This follows using sequence \eqref{five}
and the previous results.
\end{proof}

We may determine exactly the value of the right hand side of both equations
using induction and the previous lemmas. 
%

\section{$G_1B$ extensions between the simples}
We now use the $G_1$ results to classify the $G_1B$ and the $G_1T$ 
extensions between the simple $G_1B$ modules.
We use the following.
\begin{propn}
Let $\lambda$, $\mu \in X$.\\ 
\noindent(i)\ If $\mu'' - \lambda '' \in X^+$, then
$$\Ext^1_{G_1B}(\hat{L}(\lambda),\hat{L}(\mu))
\cong 
\Ext^1_{G}(L(\lambda'),L(\mu')\otimes \nabla(\mu'' -\lambda'')\frob)
$$

\noindent(ii)\ Suppose $\mu'' - \lambda '' \not\in X^+$. If $\lambda' =
\mu'$ and there exists $\alpha \in S$ and $i\in \N$ with
$\mu''-\lambda'' = -l^i\alpha$, then
$\Ext^1_{G_1B}(\hat{L}(\lambda),\hat{L}(\mu))
\cong k$.
Otherwise 
$\Ext^1_{G_1B}(\hat{L}(\lambda),\hat{L}(\mu))=0$.
\end{propn}
\begin{proof}
The proof of this proposition follows exactly as in the classical case
{\cite[proposition 9.21]{jantz2}}
\end{proof}

\begin{lem}\label{lem:lthree}
Let $\eta \in X_1$, $\mu \in X^+$. Then
$\Ext^1_G(L(\eta), 
L(\eta) \otimes \nabla(\mu)\frob)
\cong 0$.
\end{lem}
\begin{proof}
We apply the Lyndon-Hochschild-Serre five term exact sequence to
this group. Since $\Ext^1_{G_1}(L(\eta), L(\eta)) \cong 0$
we have 
$\Ext^1_G(L(\eta), 
L(\eta) \otimes \nabla(\mu)\frob)
\cong
\Ext^1_{G/G_1}(k, \nabla(\mu)\frob)
\cong
\Ext^1_{G}(k, \nabla(\mu))
\cong 0$.
\end{proof}

\begin{lem}\label{lem:little}
Let $\eta$, $\zeta \in X_1$, with $\eta\ne\zeta$ and $\mu \in X^+$. Then
$\Ext^1_G(L(\eta), 
L(\zeta) \otimes \nabla(\mu)\frob)
\cong \Hom_{G/G_1}(k, \Ext^1_{G_1}(L(\eta), L(\zeta))\otimes
\nabla(\mu)\frob)$.
\end{lem}
\begin{proof}
We apply the Lyndon-Hochschild-Serre five term exact sequence to
this group. Since $\Hom^1_{G_1}(L(\eta), L(\zeta)) \cong 0$
we have 
$\Ext^1_G(L(\eta), 
L(\eta) \otimes \nabla(\mu)\frob)
\cong \Hom_{G/G_1}(k, \Ext^1_{G_1}(L(\eta), L(\zeta))\otimes
\nabla(\mu)\frob)$.
\end{proof}

We now apply these results to our case with $G = q$-$\GL_3(k)$ or
$G=\GL_3(k)$.
We wish to determine all the extensions between the simples that
appear in a $\hat{Z}(\mu)$.
Note that the tables below will not be symmetric, we do not have
$\Ext^i_{G_1B}(\hat{L}(\mu), \hat{L}(\lambda)) \cong
\Ext^i_{G_1B}(\hat{L}(\lambda), \hat{L}(\mu))$ in general.

\begin{thm} \label{thm:g1bext}
(i)\ 
Let $(r,s) \in X_1$ with $r+s=l-2$.
If $\mu= l(a,b) +(l-1,r)$ 
then
$\Ext^1_{G_1B}(\hat{L}(\lambda), \hat{L}(\eta))$ with
$\hat{L}(\lambda)$ and $\hat{L}(\eta)$ composition factors of
$\hat{Z}(\mu)$
is given by the following table. 
 \\\medskip
\noindent\begin{tabular}{c|cccc} 
$\lambda$ $\backslash$ $\eta$
&$\mu$& $l(a-1,b)+(r,s)$ & $l(a+1,b-1)+(r,s)$& $l(a,b-1)+(s,l-1)$ \\\hline 
$\mu$ &$0$&$0$&$0$&$0$\\
$l(a-1,b)+(r,s)$ &$k$&$0$&$0$&$0$\\
$l(a+1,b-1)+(r,s)$&$0$&$k$&$0$&$0$\\
$l(a,b-1)+(s,l-1)$ &$0$&$0$&$k$&$0$\\
\end{tabular}\\
\medskip

\noindent(ii)\ 
Let $(r,s) \in X_1$ with $r+s=l-2$.
If $\mu= l(a,b) +(s,l-1)$ 
then
$\Ext^1_{G_1B}(\hat{L}(\lambda), \hat{L}(\eta))$ with
$\hat{L}(\lambda)$ and $\hat{L}(\eta)$ composition factors of
$\hat{Z}(\mu)$
is given by the following table. 
\\\medskip
\noindent\begin{tabular}{c|cccc} 
$\lambda$ $\backslash$ $\eta$
&$\mu$& $l(a,b-1)+(r,s)$ & $l(a-1,b+1)+(r,s)$& $l(a-1,b)+(l-1,r)$ \\\hline 
$\mu$ &$0$&$0$&$0$&$0$\\
$l(a-1,b)+(r,s)$ &$k$&$0$&$0$&$0$\\
$l(a+1,b-1)+(r,s)$&$0$&$k$&$0$&$0$\\
$l(a,b-1)+(s,l-1)$ &$0$&$0$&$k$&$0$\\
\end{tabular}\\
\medskip

\noindent(iii)\ 
Let $(r,s) \in X_1$ with $r+s=l-2$.
If $\mu=l(a,b) +(r,s)$ 
then
$\Ext^1_{G_1B}(\hat{L}(\lambda), \hat{L}(\eta))$ with
$\hat{L}(\lambda)$ and $\hat{L}(\eta)$ composition factors of
$\hat{Z}(\mu)$
is given by the following table. 
\\\medskip
\noindent\begin{tabular}{c|cccc} 
$\lambda$ $\backslash$ $\eta$
&$\mu$& $l(a,b-1)+(l-1,r)$ & $l(a-1,b)+(s,l-1)$& $l(a-1,b-1)+(r,s)$ \\\hline 
$\mu$ &$0$&$0$&$0$&$0$\\
$l(a-1,b)+(r,s)$ &$k$&$0$&$0$&$0$\\
$l(a+1,b-1)+(r,s)$&$k$&$0$&$0$&$0$\\
$l(a,b-1)+(s,l-1)$ &$0$&$k$&$k$&$0$\\
\end{tabular}\\
\medskip

\noindent(iv)\ 
For $(r,s) \in C$ and if $\mu = l(a,b) +(r,s)$ 
then
$\Ext^1_{G_1B}(\hat{L}(\lambda), \hat{L}(\eta))$ with
$\hat{L}(\lambda)$ and $\hat{L}(\eta)$ composition factors of
$\hat{Z}(\mu)$
is given by the following table. 
\\ \medskip
\noindent\begin{tabular}{c|ccccccccc} 
$\lambda$ $\backslash$ $\eta$
&$\mu$ & $\mu_2$ & $\mu_3$ & $\mu_4$ & $\mu_5$  
     & $\mu_6$ & $\mu_7$  & $\mu_8$ & $\mu_9$ \\ \hline
$\mu$   &$0$ & $0$ &$0$ &$0$ &$0$ &$0$ &$0$ &$0$ &$0$  \\
$\mu_2$ &$k$ & $0$ &$0$ &$0$ &$0$ &$k$ &$0$ &$0$ &$0$  \\
$\mu_3$ &$0$ & $k$ &$0$ &$0$ &$0$ &$0$ &$0$ &$0$ &$0$  \\
$\mu_4$ &$k$ & $0$ &$0$ &$0$ &$0$ &$0$ &$0$ &$k$ &$0$  \\
$\mu_5$ &$0$ & $0$ &$0$ &$k$ &$0$ &$0$ &$0$ &$0$ &$0$  \\
$\mu_6$ &$0$ & $k$ &$0$ &$0$ &$k$ &$0$ &$0$ &$0$ &$0$  \\
$\mu_7$ &$0$ & $k$ &$0$ &$k$ &$0$ &$0$ &$0$ &$0$ &$k$  \\
$\mu_8$ &$0$ & $0$ &$k$ &$k$ &$0$ &$0$ &$0$ &$0$ &$0$  \\
$\mu_9$ &$0$ & $0$ &$0$ &$0$ &$0$ &$k$ &$k$ &$k$ &$0$  
\end{tabular}\\
\medskip

\noindent(v)\ 
For $(r,s) \in C$ and if $\mu = l(a,b) +(l-s-2,l-r-2)$ 
then
$\Ext^1_{G_1B}(\hat{L}(\lambda), \hat{L}(\eta))$ with
$\hat{L}(\lambda)$ and $\hat{L}(\eta)$ composition factors of
$\hat{Z}(\mu)$
is given by the following table. 
\\
\medskip
\noindent\begin{tabular}{c|ccccccccc} 
$\lambda$ $\backslash$ $\eta$
&$\mu_1$ & $\mu_2$ & $\mu_3$ & $\mu$ & $\mu_5$  
     & $\mu_6$ & $\mu_7$  & $\mu_8$ & $\mu_9$ \\ \hline
$\mu_1$ &$0$ & $0$ &$0$ &$0$ &$0$ &$0$ &$k$ &$0$ &$0$  \\
$\mu_2$ &$k$ & $0$ &$0$ &$0$ &$k$ &$0$ &$0$ &$k$ &$0$  \\
$\mu_3$ &$0$ & $k$ &$0$ &$0$ &$0$ &$0$ &$0$ &$0$ &$k$  \\
$\mu$   &$k$ & $0$ &$0$ &$0$ &$0$ &$0$ &$0$ &$k$ &$0$  \\
$\mu_5$ &$0$ & $k$ &$0$ &$k$ &$0$ &$0$ &$0$ &$0$ &$0$  \\
$\mu_6$ &$0$ & $0$ &$0$ &$0$ &$k$ &$0$ &$0$ &$0$ &$0$  \\
$\mu_7$ &$0$ & $0$ &$0$ &$k$ &$0$ &$0$ &$0$ &$0$ &$k$  \\
$\mu_8$ &$0$ & $0$ &$0$ &$k$ &$0$ &$0$ &$0$ &$0$ &$0$  \\
$\mu_9$ &$0$ & $0$ &$0$ &$0$ &$0$ &$k$ &$k$ &$k$ &$0$  
\end{tabular}\\
\end{thm}
\begin{proof}
Most of the $\Ext$ groups above can be computed in a straight-forward
manner using the previous results.

We do sometimes need to argue as in the following case for $l=3$.

Suppose we are considering case (iv).
If $\lambda = \mu_9 = l(a-1,b-1) +(l-s-2,l-r-2)$ 
then $\mu'' - (a-1,b-1) = (1,1)$
and so
$\Ext^1_{G_1B}(\hat{L}(\lambda), \hat{L}(\mu))\cong 
\Ext^1_G(L(l-s-2,l-r-2), L(r,s)\otimes \nabla(1,1)\frob) \cong 
0$  
using lemma \ref{smallextdownalc} if $l \ge 4$.

If $l=3$ we then use the Lyndon-Hochschild-Serre five term exact
sequence.
We get
\begin{align*}
0 &\to \Ext^1_{G/G_1}(k,\Hom_{G_1}(L(l-s-2,l-r-2),L(r,s))\otimes
\bnabla(1,1)\frob)\\ &\to
\Ext^1_G(L(l-s-2,l-r-2), L(r,s)\otimes \bnabla(1,1)\frob)\\ &\to
\Hom_{G/G_1}(k,\Ext^1_{G_1}(L(l-s-2,l-r-2),L(r,s))\otimes
\bnabla(1,1)\frob)\\ &\to
\Ext^2_{G/G_1}(k,\Hom_{G_1}(L(l-s-2,l-r-2),L(r,s))\otimes
\bnabla(1,1)\frob)\\ &\to
\Ext^2_G(L(l-s-2,l-r-2), L(r,s)\otimes \bnabla(1,1)\frob) 
\end{align*}
Since $\Hom_{G_1}(L(l-s-2,l-r-2),L(r,s))$ is zero we have
using theorem \ref{thm:yehia}
\begin{align*}
  \Ext^1_G(&L(l-r-2,l-s-2), L(r,s)\otimes \bnabla(1,1)\frob)\\
&\cong
\Hom_{G/G_1}(k,\Ext^1_{G_1}(L(l-r-2,l-s-2),L(r,s))\otimes
\nabla(1,1)\frob)\\ 
&\cong
\Hom_{G/G_1}(k,
\bnabla(1,1)\frob \oplus \bnabla(0,1)\frob \otimes \bnabla(1,1)\frob
\oplus \bnabla(1,0)\frob \otimes \bnabla(1,1)\frob
)\\ 
&\cong
\Hom_{\SL_3}(k,
\bnabla(1,1) \oplus \bnabla(0,1) \otimes \bnabla(1,1)
\oplus \bnabla(1,0) \otimes \bnabla(1,1)
)\\ 
&\cong
\Hom_{\SL_3}(k,
\bnabla(1,1))\oplus 
\Hom_{G}(\bnabla(1,0), \bnabla(1,1))
\oplus \Hom_{G}(\bnabla(0,1), \bnabla(1,1))
\\ 
&\cong 0
\qedhere
\end{align*}

%
%
%
%
%
\end{proof}

\section{The composition series of induced modules for $G_1B$.}

Before deducing the $G_1B$ structure of the $\hat{Z}(\mu)$'s we need some
more propositions.

\begin{propn}\label{propn:ztrans}
Suppose $\lambda \in C$ and $\mu  \in \bar{C}$ lies on a
wall. Suppose also that $s$ is a simple reflection that fixes $\mu$.
and that $w \cdot \lambda < ws\cdot \lambda$.
We have the following properties.
\begin{enumerate}
\item[(i)]
$T_{\lambda}^{\mu} L(w \cdot \lambda) \cong L(w \cdot \mu)$
and 
$T_{\lambda}^{\mu} L(ws \cdot \lambda) \cong 0$
\item[(ii)]
$\hat{T}_\lambda^\mu \hat{Z}(w \cdot \lambda)
\cong \hat{T}_\lambda^\mu \hat{Z}(ws \cdot \lambda)
\cong \hat{Z}(w \cdot \mu)
$
\item[(iii)]
We have a short exact sequence
$$ 0 \to \hat{Z}(w\cdot\lambda) \to 
\hat{T}_\mu^\lambda \hat{Z}(w \cdot \mu)
\to \hat{Z}(ws\cdot\lambda) \to 0$$
The socle of $\hat{T}_\mu^\lambda \hat{Z}(w \cdot \mu)$ is
$\hat{L}(w\cdot\lambda)$.
\end{enumerate}

\end{propn}
This is the quantum version of \cite[II 9.22 (4), (2), (3)]{jantz2}
and may be proved as in the classical case using the results of
\cite{andpolwen} and \cite{donkbk}.

\begin{propn}\label{propn:ztranshom}
Let $\lambda$, $\mu$ , $w$ and $s$ be as in the previous proposition.
We have
$\Hom_{G_1B}(\hat{Z}(ws\cdot\lambda), \hat{Z}(w\cdot\lambda)) \cong 
\Hom_{G_1B}(\hat{Z}(w\cdot\mu), \hat{Z}(w\cdot\mu)) \cong 
k
$ 
\end{propn}
\begin{proof}
Firstly,
we have
$\Hom_{G_1B}(\hat{Z}(w\cdot\mu), \hat{Z}(w\cdot\mu)) \cong 
\Hom_{B}(\hat{Z}(w\cdot\mu), k_{w\cdot\mu}) 
$
by Frobenious reciprocity. The latter group is at most one
dimensional, as the dimension of the $w \cdot \mu$ weight space in
$\hat{Z}(w\cdot\mu)$ is one.
On the other hand 
$\Hom_{G_1B}(\hat{Z}(w\cdot\mu), \hat{Z}(w\cdot\mu))$ is certainly
non-zero. Thus there is unique homomorphism (upto scalars), the identity
homomorphism.

We may now argue as in the proof of \cite[II, proposition 7.19]{jantz2}
to show that the map $\phi$ in the following long exact sequence is
zero,
$$ 
0 \to \Hom_{G_1B}(\hat{Z}(ws\cdot\lambda), \hat{Z}(w\cdot\lambda)) \to 
\Hom_{G_1B}(\hat{Z}(ws\cdot\lambda), \hat{T}_\mu^\lambda \hat{Z}(w \cdot \mu))
\stackrel{\phi}\to \Hom_{G_1B}(\hat{Z}(ws\cdot\lambda), \hat{Z}(ws\cdot\lambda))$$
and we thus get the isomorphism as claimed.
%
\end{proof}

We may now prove the following theorem,
We use the following various facts about $\hat{Z}(\lambda)$ for
$\lambda \in X^+$:
\begin{enumerate}
\item[(i)]
$\hat{Z}(\lambda)$ has simple $G_1B$ socle $\hat{L}(\lambda)$
(see \cite[II, 9.6 (1)]{jantz2} and \cite[3.1 (13) (i)]{donkbk})
\item[(ii)]
$\hat{Z}(\lambda)$ has simple $G_1B$ head $\hat{L}(2(l-1)\rho
  -\lambda)^* \cong \hat{L}(2(l-1)\rho +w_0\lambda + l(w_0\lambda''
  -\lambda'))
$
(see \cite[II, 9.6 (2)]{jantz2} and \cite[3.1 (22)]{donkbk})
\item[(iii)]
$\hat{Z}(\lambda)^*
\cong 
\hat{Z}(2(l-1)\rho -\lambda)
$
(see  \cite[II, 9.2 (2)]{jantz2} and \cite[3.1 (21)]{donkbk})
\item[(iv)]
$\hat{Z}(\lambda + l\mu)
\cong 
\hat{Z}(\lambda)\otimes k_{l\mu}
$ 
(see \cite[II, 9.2 (5)]{jantz2}, also follows in the quantum case using
the tensor identity).
\end{enumerate}
Strictly speaking the results in the quantum case using \cite{donkbk}
are only $G_1T$ results. But the above properties clearly lift to
$G_1B$.

Item (iii) above implies that the submodule structure of $\hat{Z}(\lambda)$ for
$\lambda$ in a down alcove and the structure of $\hat{Z}(\mu)$ for
$\mu$ in an up alcove are inversions of each other.
Item (iv) above implies that the structure for a weight of a
particular $G_1$ type is always the same.


\begin{thm}\label{thm:hatzs}
The submodule structure of the $\hat{Z}(\lambda)$ for $\lambda \in
X^+$ is as follows.
\begin{enumerate}
\item[(i)] Suppose 
$\lambda = l(a,b) + (l-1,l-1)$ with $(a,b) \in X^+$. Then
$$\hat{Z}(\lambda)= \hat{L}(\lambda).$$
\item[(ii)] Suppose
$\lambda = l(a,b) + (l-1,r)$ with $(a,b) \in X^+$ and $(l-1,r)\in
X_1$.
If we set $s=l-r-2$ then 
the module $\hat{Z}(\lambda)$ has filtration 
$$\xymatrix@R=10pt
{
{\hat{L}(l(a,b-1)+(s,l-1))} \ar@{-}[d] \\
{\hat{L}(l(a+1,b-1)+(r,s))} \ar@{-}[d] \\
{\hat{L}(l(a-1,b)+ (r,s))} \ar@{-}[d] \\
{\hat{L}(l(a,b)+ (l-1,r)).}}
$$ 
\item[(iii)]\  Suppose
$\lambda = l(a,b) + (s,l-1)$ with $(a,b) \in X^+$ and $(s,l-1)\in
X_1$.
If we set $r=l-s-2$ then the module
$\hat{Z}(\lambda)$ has filtration 
$$\xymatrix@R=10pt
{
{\hat{L}(l(a-1,b)+ (l-1,r))} \ar@{-}[d] \\
{\hat{L}(l(a-1,b+1)+ (r,s))} \ar@{-}[d] \\
{\hat{L}(l(a,b-1)+ (r,s))} \ar@{-}[d] \\
{\hat{L}(l(a,b)+ (s,l-1)).}}
$$ 
\item[(iv)]\ Suppose
$\lambda = l(a,b) + (r,s)$ with $(a,b) \in X^+$, $(r,s)\in X_1$ and
$r+s=l-2$. Then the module
$\hat{Z}(\lambda)$ has filtration
$$\xymatrix@R=10pt@C=-20pt
{
& {\hat{L}(l(a-1,b-1)+ (r,s))} \ar@{-}[dl] \ar@{-}[dr]
\\
{\hat{L}(l(a,b-1)+ (l-1,r))} \ar@{-}[dr] &
 &{\hat{L}(l(a-1,b)+ (s,l-1))} \ar@{-}[dl] \\
& {\hat{L}(l(a,b)+ (r,s)).}}
$$
\item[(v)]\ Suppose
$\lambda = l(a,b) + (r,s)$ with $(a,b) \in X^+$, 
and $(r,s)\in C$. 
We let $\mu_1$ upto $\mu_9$ be as before,
 depicted in Figure~1~(a), where the number corresponds to
the subscript of $\mu$.

Then $\hat{Z}(\lambda)$ has filtration 
$$\xymatrix@R=10pt@C=5pt
{
& &{\hat{L}(\mu_9)} \ar@{-}[dll] 
\ar@{-}[drr] \ar@{-}[d] \\
{{\hat{L}(\mu_6)}} \ar@{-}[ddrrr]   \ar@{-}[d]   
& & {\hat{L}(\mu_7)} \ar@{-}[ddr] \ar@{-}[ddl] 
& & {\hat{L}(\mu_8)} \ar@{-}[ddlll] \ar@{-}[d] \\
{\hat{L}(\mu_5)} \ar@{-}[dr]
& & & &{\hat{L}(\mu_3)} \ar@{-}[dl] \\
& {\hat{L}(\mu_4)} \ar@{-}[dr]
& &{\hat{L}(\mu_2)} \ar@{-}[dl] \\
& &{\hat{L}(\mu_1).}}
$$
\item[(vi)]\ Suppose
$\lambda = l(a,b) + (l-s-2,l-r-2)$ with $(a,b) \in X^+$, 
and $(r,s)\in C$. 
We let $\mu_1$ upto $\mu_9$ be as before,
 depicted in Figure~1~(b), where the number corresponds to
the subscript of $\mu$.
Then $\hat{Z}(\lambda)$ has filtration 
$$\xymatrix@R=10pt@C=5pt
{
& &{\hat{L}(\mu_3)} \ar@{-}[dl] \ar@{-}[dr]  \\
& {\hat{L}(\mu_2)} \ar@{-}[ddr]\ar@{-}[ddrrr]\ar@{-}[dl] 
& &{\hat{L}(\mu_9)} \ar@{-}[dr] \ar@{-}[ddl]\ar@{-}[ddlll] \\
{\hat{L}(\mu_1)} \ar@{-}[d]
& & & &{\hat{L}(\mu_6)} \ar@{-}[d] \\
{{\hat{L}(\mu_7)}} \ar@{-}[drr]   
& & {\hat{L}(\mu_8)} \ar@{-}[d] 
& & {\hat{L}(\mu_5)} \ar@{-}[dll] \\
& &{\hat{L}(\mu_4).}}
$$
\end{enumerate}
\end{thm}
\begin{proof}
The structures for (i)-(iv) are the only possible ones using the fact
that $\hat{Z}(\lambda)$ has simple head and socle as described above
and the possible extensions that exist between the composition
factors.

Cases (v) and (vi).
The structure depicted has all the possible extensions drawn in. We
need to prove that all these extensions do actually appear.
The simples must be in the layers as described, for otherwise it would
contradict the $\hat{Z}(\lambda)$ having simple socle
$\hat{L}(\lambda)$ and simple head $L(\mu_9)$ ($L(\mu_3)$) if
$\lambda$ is a down (up) alcove respectively.

For instance, in case (v) we must have a uniserial subquotient of
$\hat{L}(\mu_4)$, $\hat{L}(\mu_5)$ and $\hat{L}(\mu_6)$, since
$\hat{L(\mu_5)}$ can only extend one simple below it (namely
$\hat{L}(\mu_4)$) and one simple above it, (namely $\hat{L}(\mu_6)$). Otherwise
$\hat{L}(\mu_5)$ would either be in the head or socle of
$\hat{Z}(\lambda)$.

So for case (v) we can deduce the following structure so far:
$$\xymatrix@R=10pt@C=5pt
{
& &{\hat{L}(\mu_9)} \ar@{-}[dll] 
\ar@{-}[drr] \ar@{-}[d] \\
{{\hat{L}(\mu_6)}} 
 \ar@{-}[d]   
& & {\hat{L}(\mu_7)} 
& & {\hat{L}(\mu_8)} 
\ar@{-}[d] \\
{\hat{L}(\mu_5)} \ar@{-}[dr]
& & & &{\hat{L}(\mu_3)} \ar@{-}[dl] \\
& {\hat{L}(\mu_4)} \ar@{-}[dr]
& &{\hat{L}(\mu_2)} \ar@{-}[dl] \\
& &{\hat{L}(\mu_1).}}
$$
We get a similar picture (only inverted) for case (vi).
$$\xymatrix@R=10pt@C=5pt
{
& &{\hat{L}(\mu_3)} \ar@{-}[dl] \ar@{-}[dr]  \\
& {\hat{L}(\mu_2)} 
\ar@{-}[dl] 
& &{\hat{L}(\mu_9)} \ar@{-}[dr] 
\\
{\hat{L}(\mu_1)} \ar@{-}[d]
& & & &{\hat{L}(\mu_6)} \ar@{-}[d] \\
{{\hat{L}(\mu_7)}} \ar@{-}[drr]   
& & {\hat{L}(\mu_8)} \ar@{-}[d] 
& & {\hat{L}(\mu_5)} \ar@{-}[dll] \\
& &{\hat{L}(\mu_4).}}
$$

Consider the structure for case (v) so far.
The $\hat{L}(\mu_7)$ must extend at least one of $\hat{L}(\mu_4)$ or
$\hat{L}(\mu_2)$. Suppose that it extends  $\hat{L}(\mu_4)$.
Now the existence of a homomorphism from $\hat{Z}(\mu_1)$
to $\hat{Z}(\mu_4)$ (using proposition \ref{propn:ztranshom}) implies that 
there is an extension of $\hat{L}(\mu_7)$ by $\hat{L}(\mu_9)$ in
$\hat{Z}(\mu_4)$, as the image of the homorphism must contain at least
$\hat{L}(\mu_4)$, $\hat{L}(\mu_5)$, 
$\hat{L}(\mu_6)$, $\hat{L}(\mu_7)$ and
$\hat{L}(\mu_9)$, and it has simple head $\hat{L}(\mu_9)$.

Now consider the module $\hat{Z}(\eta)$ defined to be 
$\hat{Z}(\mu_4)^* \otimes k_{l(2a-1,2b-1)}$. The weight $\eta$ is in
the same (down) alcove as the $\mu_8$ from $\hat{Z}(\mu_4)$.
We now consider the dual of the extension of 
$\hat{L}(\mu_7)$ by $\hat{L}(\mu_9)$  and tensor it by 
$ k_{l(2a-1,2b-1)}$. This extension then appears in
$\hat{Z}(\eta)$ and working out what the duals of the simples are
gives us an extension of $\hat{L}(\eta_4)$ by
$\hat{L}(\eta_8)$. Translation principle then tells us that our
original $\hat{Z}(\mu_1)$ has an extension of  $\hat{L}(\mu_4)$ by
$\hat{L}(\mu_8)$. 

Considering the homomorphism from $\hat{Z}(\mu_1)$
to $\hat{Z}(\mu_4)$ again implies that 
there is an extension of $\hat{L}(\mu_8)$ by $\hat{L}(\mu_9)$ 
in $\hat{Z}(\mu_4)$. 

So we now have for case (v) (assuming that $\hat{L}(\mu_4)$ extends
$\hat{L}(\mu_7)$)
$$\xymatrix@R=10pt@C=5pt
{
& &{\hat{L}(\mu_9)} \ar@{-}[dll] 
\ar@{-}[drr] \ar@{-}[d] \\
{{\hat{L}(\mu_6)}} 
 \ar@{-}[d]   
& & {\hat{L}(\mu_7)} \ar@{-}[ddl] 
& & {\hat{L}(\mu_8)} \ar@{-}[ddlll] 
\ar@{-}[d] \\
{\hat{L}(\mu_5)} \ar@{-}[dr]
& & & &{\hat{L}(\mu_3)} \ar@{-}[dl] \\
& {\hat{L}(\mu_4)} \ar@{-}[dr]
& &{\hat{L}(\mu_2)} \ar@{-}[dl] \\
& &{\hat{L}(\mu_1).}}
$$
For case (vi) we get:
$$\xymatrix@R=10pt@C=5pt
{
& &{\hat{L}(\mu_3)} \ar@{-}[dl] \ar@{-}[dr]  \\
& {\hat{L}(\mu_2)} 
\ar@{-}[dl] 
& &{\hat{L}(\mu_9)} \ar@{-}[dr] \ar@{-}[ddl]\ar@{-}[ddlll] 
\\
{\hat{L}(\mu_1)} \ar@{-}[d]
& & & &{\hat{L}(\mu_6)} \ar@{-}[d] \\
{{\hat{L}(\mu_7)}} \ar@{-}[drr]   
& & {\hat{L}(\mu_8)} \ar@{-}[d] 
& & {\hat{L}(\mu_5)} \ar@{-}[dll] \\
& &{\hat{L}(\mu_4).}}
$$
Now the image of the
homomorphism from $\hat{Z}(\mu_4)$ to $\hat{Z}(\mu_8)$ (which exists
using proposition \ref{propn:ztranshom}) contains an extension of 
$\hat{L}(\mu_9)$ and $\hat{L}(\mu_3)$.  Thus there is also an
an extension of $\hat{L}(\mu_2)$ and $\hat{L}(\mu_7)$ 
in the original $\hat{Z}(\mu_1)$ for case (v). 

Repeating the above argument with $\mu_2$ in place of $\mu_4$ thus
gives us the result.
\end{proof}

\section{The good $l$-filtrations of the induced modules for $G$}
\begin{thm}\label{thm:lfiltra}
Each $\nabla(\lambda)$ has a $l$-filtration.
This filtration takes the following form:
\begin{enumerate}
\item[(i)] Suppose 
$\lambda = l(a,b) + (l-1,l-1)$ with $(a,b) \in X^+$. Then
$$\nabla(\lambda)= \bnabla(a,b)\frob\otimes L(l-1,l-1).$$
\item[(ii)] Suppose
$\lambda = l(a,b) + (l-1,r)$ with $(a,b) \in X^+$ and $(l-1,r)\in
X_1$.
If we set $s=l-r-2$ then for $a \equiv -1 \pmod l$,
the module $\nabla(\lambda)$ has filtration 
$$\xymatrix@R=10pt
{
{\bnabla(a,b-1)\frob\otimes L(s,l-1)} \ar@{-}[d] \\
{\bnabla(a+1,b-1)\frob\otimes L(r,s)} \ar@{-}[d] \\
{\bnabla(a-1,b)\frob\otimes L(r,s)} \ar@{-}[d] \\
{\bnabla(a,b)\frob\otimes L(l-1,r)}}
$$ 
while for  $a \not\equiv -1 \pmod l$,
$\nabla(\lambda)$ has filtration 
$$\xymatrix@R=10pt@C=-20pt
{
& {\bnabla(a,b-1)\frob\otimes L(s,l-1)} \ar@{-}[dl] \ar@{-}[dr]
\\
{\bnabla(a+1,b-1)\frob\otimes L(r,s)} \ar@{-}[dr] &
 &{\bnabla(a-1,b)\frob\otimes L(r,s)} \ar@{-}[dl] \\
& {\bnabla(a,b)\frob\otimes L(l-1,r).}}
$$ 
\item[(iii)]\  Suppose
$\lambda = l(a,b) + (s,l-1)$ with $(a,b) \in X^+$ and $(s,l-1)\in
X_1$.
If we set $r=l-s-2$ then for $b \equiv -1 \pmod l$, the module
$\nabla(\lambda)$ has filtration 
$$\xymatrix@R=10pt
{
{\bnabla(a-1,b)\frob\otimes L(l-1,r)} \ar@{-}[d] \\
{\bnabla(a-1,b+1)\frob\otimes L(r,s)} \ar@{-}[d] \\
{\bnabla(a,b-1)\frob\otimes L(r,s)} \ar@{-}[d] \\
{\bnabla(a,b)\frob\otimes L(s,l-1)}}
$$ 
while for  $b \not\equiv -1 \pmod l$,
$\nabla(\lambda)$ has filtration 
$$\xymatrix@R=10pt@C=-20pt
{
& {\bnabla(a-1,b)\frob\otimes L(l-1,r)} \ar@{-}[dl] \ar@{-}[dr]
\\
{\bnabla(a-1,b+1)\frob\otimes L(r,s)} \ar@{-}[dr] &
 &{\bnabla(a,b-1)\frob\otimes L(r,s)} \ar@{-}[dl] \\
& {\bnabla(a,b)\frob\otimes L(s,l-1).}}
$$ 
\item[(iv)]\ Suppose
$\lambda = l(a,b) + (r,s)$ with $(a,b) \in X^+$, $(r,s)\in X_1$ and
$r+s=l-2$. Then the module
$\nabla(\lambda)$ has filtration
$$\xymatrix@R=10pt@C=-20pt
{
& {\bnabla(a-1,b-1)\frob\otimes L(r,s)} \ar@{-}[dl] \ar@{-}[dr]
\\
{\bnabla(a,b-1)\frob\otimes L(l-1,r)} \ar@{-}[dr] &
 &{\bnabla(a-1,b)\frob\otimes L(s,l-1)} \ar@{-}[dl] \\
& {\bnabla(a,b)\frob\otimes L(r,s).}}
$$
\item[(v)]\ Suppose 
$\lambda = l(a,0) + (r,s)$ with $(a,0) \in X^+$, $a\ge1$ and $(r,s)\in
  C$ 
then the module $\nabla(\lambda)$ has filtration
$$\xymatrix@R=10pt{
{\bnabla(a-2,0)\frob\otimes L(s,l-r-s-3)} \ar@{-}[d] \\
{\bnabla(a-1,0)\frob\otimes L(l-r-2,r+s+1)} \ar@{-}[d] \\
{\bnabla(a,0)\frob\otimes L(r,s).}}
$$
\item[(vi)]\ Suppose
$\lambda = l(0,b) + (r,s)$ with $(0,b) \in X^+$, $b\ge1$ and $(r,s)\in
  C$.
Then the module $\nabla(\lambda)$ has filtration
$$\xymatrix@R=10pt{
{\bnabla(0,b-2)\frob\otimes L(l-r-s-3,r)} \ar@{-}[d] \\
{\bnabla(0,b-1)\frob\otimes L(r+s+1,l-s-2)} \ar@{-}[d] \\
{\bnabla(0,b)\frob\otimes L(r,s).}}
$$
\item[(vii)]\ Suppose
$\lambda = l(a,b) + (r,s)$ with $(a,b) \in X^+$, $a$ and $b\ge 1$, 
and $(r,s)\in C$. 
We let $\mu_1$ upto $\mu_9$ be as before,
 depicted in Figure~1~(a), where the number corresponds to
the subscript of $\mu$.

Then for $a$ and $b \equiv 0
\pmod l$, $\nabla(\lambda)$ has filtration 
$$\xymatrix@R=10pt@C=5pt
{
& &{\nabla_l(\mu_9)} \ar@{-}[dll] 
\ar@{-}[drr] \ar@{-}[d] \\
{{\nabla_l(\mu_6)}} \ar@{-}[ddrrr]   \ar@{-}[d]   
& & {\nabla_l(\mu_7)} \ar@{-}[ddr] \ar@{-}[ddl] 
& & {\nabla_l(\mu_8)} \ar@{-}[ddlll] \ar@{-}[d] \\
{\nabla_l(\mu_5)} \ar@{-}[dr]
& & & &{\nabla_l(\mu_3)} \ar@{-}[dl] \\
& {\nabla_l(\mu_4)} \ar@{-}[dr]
& &{\nabla_l(\mu_2)} \ar@{-}[dl] \\
& &{\nabla_l(\mu_1).}}
$$
For $a  \not\equiv 0 \pmod l$ there is no extension of
$\nabla_l(\mu_5)$ by $\nabla_l(\mu_6)$.
For $b  \not\equiv 0 \pmod l$ there is no extension of
$\nabla_l(\mu_3)$ by $\nabla_l(\mu_8)$.
So for $a$ and $b \not\equiv 0 \pmod l$ we have:
$$\xymatrix@R=10pt@C=5pt
{
& & &{\nabla_l(\mu_9)} \ar@{-}[dlll] \ar@{-}[drrr] \ar@{-}[dll] 
                \ar@{-}[drr] \ar@{-}[d] \\
{\nabla_l(\mu_5)} \ar@{-}[drr]
& {\ \nabla_l(\mu_6)\!\!\!} \ar@{-}[drrr] \ar@{-}[dr]    
& & {\nabla_l(\mu_7)} \ar@{-}[dr] \ar@{-}[dl] 
& & {\!\!\!\nabla_l(\mu_8)\ } \ar@{-}[dlll] \ar@{-}[dl] 
& {\nabla_l(\mu_3)} \ar@{-}[dll] \\
& &{{\!\nabla_l(\mu_4)\!}} \ar@{-}[dr]
& & {{\!\nabla_l(\mu_2)\!}} \ar@{-}[dl] \\
& & &{\nabla_l(\mu_1)}}
$$
and similarly for the other cases for $a$ and $b$.
\item[(viii)]\ Suppose
$\lambda = l(a,b) + (l-s-2,l-r-2)$ with $(a,b) \in X^+$, 
and $(r,s)\in C$. 
We let $\mu_1$ upto $\mu_9$ be as before,
 depicted in Figure~1~(b), where the number corresponds to
the subscript of $\mu$.
Then for $a$ and $b \equiv -1
\pmod l$, $\nabla(\lambda)$ has filtration 
$$\xymatrix@R=10pt@C=5pt
{
& &{\nabla_l(\mu_3)} \ar@{-}[dl] \ar@{-}[dr]  \\
& {\nabla_l(\mu_2)} \ar@{-}[ddr]\ar@{-}[ddrrr]\ar@{-}[dl] 
& &{\nabla_l(\mu_9)} \ar@{-}[dr] \ar@{-}[ddl]\ar@{-}[ddlll] \\
{\nabla_l(\mu_1)} \ar@{-}[d]
& & & &{\nabla_l(\mu_6)} \ar@{-}[d] \\
{{\nabla_l(\mu_7)}} \ar@{-}[drr]   
& & {\nabla_l(\mu_8)} \ar@{-}[d] 
& & {\nabla_l(\mu_5)} \ar@{-}[dll] \\
& &{\nabla_l(\mu_4).}}
$$
For $a \not\equiv -1\pmod l$ there is no extension of
$\nabla_l(\mu_5)$ by $\nabla_l(\mu_6)$.
For $b \not\equiv -1\pmod l$ there is no extension of
$\nabla_l(\mu_7)$ by $\nabla_l(\mu_1)$.
So for $a$ and $b \not\equiv -1\pmod l$ we have:
$$\xymatrix@R=10pt@C=5pt
{
& & &{\nabla_l(\mu_3)} \ar@{-}[dl] \ar@{-}[dr]  \\
& & {\nabla_l(\mu_2)} \ar@{-}[dr]\ar@{-}[drrr]\ar@{-}[dl]\ar@{-}[dll]  
& &{\nabla_l(\mu_9)} \ar@{-}[dr] \ar@{-}[dl]\ar@{-}[dlll]\ar@{-}[drr]  \\
{\nabla_l(\mu_1)} \ar@{-}[drrr]
&{{\ \nabla_l(\mu_7)\!\!\!}} \ar@{-}[drr]   
& &{\nabla_l(\mu_8)} \ar@{-}[d] 
& &{\!\!\!\nabla_l(\mu_5)\ } \ar@{-}[dll] 
&{\nabla_l(\mu_6)} \ar@{-}[dlll] \\
& & &{\nabla_l(\mu_4)}}
$$
and similarly for the other cases for $a$ and $b$.
\end{enumerate}
\end{thm}
\begin{proof}
This may now be proved as in the classical case \cite{parker1}.
\end{proof}

\section{Homorphisms between induced modules for $q$-$\GL_3(k)$}
We now show how to generalise the results of \cite{coxpar} to the
quantum case. As noted in that paper, there were two obstacles to 
this. The first was that we needed an $l$-filtration of the induced
modules, and the second
was that we needed a quantum version of main result of
\cite{cpmaps}. We can now prove that this result (\cite{cpmaps}) holds for
$q$-$\GL_3(k)$, but unfortunately not in general.
We will assume that $p \ne 0$. The case with $p=0$ is easier.

We define a \emph{$lp^e$-wall} for $e \in\N$ to be a wall for $X^+$ that is
fixed by a reflection of the form $s_{\beta, mlp^e}$ for some $m \in
\Z$ and $\beta \in R$.
%
%
\begin{thm}
Suppose that $\lambda, \mu \in X^+$
satisfy the following conditions:    
\begin{enumerate}
\item[(i)]{$\mu < \lambda $.  
}
\item[(ii)]{There exists some $e\in\N$ such that:} 
\begin{enumerate}
\item[(a)]{$\lambda$ and $\mu$ are mirror images in some $lp^e$-wall
$L$ and}
\item[(b)]{$L$ is the unique $lp^e$-wall between $\lambda $ and $\mu$ (possibly containing
$\lambda$ or $\mu$) parallel to $L$.}
\end{enumerate}
\end{enumerate}
Then $\Hom_G(\nabla(\lambda), \nabla(\mu)) \neq 0$.
\end{thm}
\begin{proof}
We may assume that $\lambda$ is not a Steinberg weight as then the
result follows by twisting the corresponding map for the classical
case.

Suppose $L$ is fixed by $s_{\beta, mlp^e}$ for some $m \in \N$ and
$\beta  \in \R^+$.
There are two cases to consider. 

Case (1): $\beta$ is a simple root.
In this case the theorem reduces to the analogous one for
$q$-$\GL_2(k)$ using Levi subgroups and the results of Donkin
\cite{donkbk}. See \cite[theorem 5.1 and 7.1]{coxerd}.

Case (2): $\beta=\rho$. In this case, we construct the
homomorphism directly.

We first suppose that $e=0$ and that
$\lambda$ doesn't lie in an up alcove. We then claim
that the required
map is the one obtained by inducing the map
$\hat{Z}(\lambda) \to \hd(\hat{Z}(\lambda))$
from $G_1B$ upto $G$.

We claim that the head of $\hat{Z}(\lambda)$ is
$\hat{L}(\mu)$. 
We write $\lambda = l(a,b) + (r,s)$ with  $(a,b)\in X^+$ and $(r,s)\in
X_{1}$.
Now 
\begin{align*}
 \hd(\hat{Z}(\lambda)) &= 
\hat{L}( 2(l -1)\rho -\lambda)^*\\
&=\hat{L}(l(2-a,2-b)-(r+2,s+2))^*\\
&\cong
L(-w_0(l-r-2,l-s-2))\otimes
k_{-l(1-a,1-b)}\\
&=
L(l-s-2,l-r-2))\otimes
k_{-l(1-a,1-b)}\\
&\cong
\hat{L}((l-s-2,l-r-2)+l(a-1,b-1))
\end{align*}
Also the condition on $L$, $\lambda$ and $\mu$ implies that
$m$ is the greatest integer such that
$\langle \lambda +\rho , \crho \rangle - ml$ is positive.
We thus have 
$\langle \lambda +\rho , \crho \rangle= ml + d$, where $1\le d\le
l$.
Hence
\begin{align*}
\mu
&= s_{\rho, ml} \cdot \lambda\\
&= \lambda - (\langle \lambda +\rho , \crho \rangle - ml) \rho\\
&= \lambda - d\rho.
\end{align*}
Since 
$\langle \lambda +\rho , \crho \rangle= 
l(a+b) + r+s+2$, $d$ is then $r+s+2$, as the condition that
$\lambda$ is not in an up alcove implies that $r+s+2$ is at most
$l$.
Thus $\mu = (l-s-2,l-r-2)+l(a-1,b-1)$, as required.

We note that the image of this map is $\Ind_{G_1B}^G \hat{L}(\mu) = 
\nabla_l(\mu)$.

If $e=0$ and $\lambda$ lies in an up alcove then the
required map is that of \ref{propn:ztranshom}. (It has image
the quotient module of $\nabla(\lambda)$ with an $l$-filtration by
$\nabla_l(\lambda_8)$, 
$\nabla_l(\lambda_2)$, 
$\nabla_l(\lambda_3)$ 
and 
$\nabla_l(\lambda_9)$.)

We now suppose $ e>0$.
We let $\nabla_l(\eta)$ be the $G_1$-head of $\nabla(\lambda)$. 
We know that 
$\eta = \lambda - (r+s+2)\rho$, using the same notation as in the
previous case.  Note that $\eta'=\mu'$, as they are both downward
reflections of $\lambda$.

We claim that $\eta''$, (considered as a weight for $\SL_3(k)$) is
$s_{\beta, mp^e} \cdot \mu''$.
Thus there is a Carter-Payne map from 
$$ \phi:\bnabla(\eta'') \to \bnabla(\mu'').$$
We then twist the above map:
$$Id \otimes \phi\frob : \nabla_l(\eta) \to \nabla_l(\mu).$$
This then induces the required map from $\nabla(\lambda)$ to $\nabla(\mu)$.

We now prove the claim.
Consider
$s_{\beta, mp^e} \cdot \mu''
= \mu'' - (\langle \mu'' +\rho, \crho \rangle-mp^e)\rho.$
Now the condition on $L$, $\lambda$ and $\mu$ imply that 
$\langle \mu +\rho , \crho \rangle= mlp^e - d$, where $1\le d\le
lp^e$. Thus
$\langle \mu'', \crho \rangle -  mp^e= -\frac{1}{l}(
d+\langle \mu'+\rho, \crho \rangle)$.
And so
\begin{align*}
l (s_{\beta, mp^e} \cdot \mu'') +\mu'
&= l\mu'' + (d+\langle \mu'+\rho,\crho\rangle -2l) \rho +\mu'\\
&= \mu + d \rho + (\langle \eta'+\rho,\crho\rangle -2l) \rho \\
&= \lambda - (2l-\langle \eta'+\rho,\crho\rangle) \rho \\
&= \lambda - (2l-\langle (l-s-1,l-r-1),(1,1)\rangle) \rho \\
&= \lambda - (s+r+2)\rho \\
&= \eta.
\end{align*}
Thus $\eta'' = s_{\beta, mp^e} \cdot \mu''$  as required.
\end{proof}

As a corollary we get that 
all the results of \cite{coxpar} regarding homomorphisms between 
induced modules now generalise to the quantum case if $l \ge 3$.
We just need to replace the $p^{e+1}$ walls and reflections with
$lp^e$ walls and reflections.

In particular we have
\begin{thm}
Suppose $l \ge 3$, then
all the $\Hom_G(\nabla(\lambda), \nabla(\mu))$, with $\lambda$, $\mu
\in X^+$ are at most
one-dimensional.
\end{thm}
The non-zero homomorphisms may be determined by using the appropriate
generalisations of the main theorems of \cite{coxpar}.

The characteristic zero case is easier. Here, we only get 
reflections about 
\emph{$l$-walls}, that is, a wall fixed by a reflection of the form  
$s_{\beta, ml}$ for some $m \in
\Z$ and $\beta \in R$.
In this case, the only maps are the $l$-good maps. This is because
$\nabla_l(\lambda)$ is always isomorphic to $L(\lambda)$, thus any map
between induced modules must respect the $l$-filtration. 
Hence we have the following.
\begin{thm}
Suppose $p=0$.
All the $\Hom_G(\nabla(\lambda), \nabla(\mu))$, with $\lambda$, $\mu
\in X^+$ are at most
one-dimensional.
Any non-zero map is an $l$-good map and is described by the
appropriate quantum version of 
\cite[Lemma 3.1]{coxpar}.
\end{thm}


\begin{thebibliography}{10}

\bibitem{anderpfil}
H.~H. Andersen, \emph{{$p$}--filtrations and the {S}teinberg module}, J.
  Algebra \textbf{244} (2001), 664--683.

\bibitem{andpolwen}
H.~H. Andersen, P.~Polo, and K.~X. Wen, \emph{Representations of quantum
  algebras}, Invent. Math. \textbf{104} (1991), no.~1, 1--59.

\bibitem{benson}
D.~J. Benson, \emph{{R}epresentations and {C}ohomology {I}}, Cambridge Studies
  in Advanced Mathematics, no.~30, Cambridge University Press, 1995.

\bibitem{cpmaps}
R.~Carter and M.~T.~J. Payne, \emph{On homomorphisms between weyl modules and
  specht modules}, Math. Proc. Cambridge Philos. Soc. \textbf{87} (1980),
  no.~3, 419--425.

\bibitem{coxerd}
A.~G. Cox and K.~Erdmann, \emph{On {${\mathrm{Ext}}^2$} between {W}eyl modules
  for quantum {${\mathrm{GL}}_n$}}, Math. Proc. Cambridge Philos. Soc.
  \textbf{128} (2000), 441--463.

\bibitem{coxpar}
A.~G. Cox and A.~E. Parker, \emph{Homomorphisms between {W}eyl modules for
  $\mathrm{SL}_3(k)$}, preprint, 2003.

\bibitem{dipdonk}
R.~Dipper and S.~Donkin, \emph{Quantum {$\mathrm{GL}_n$}}, Proc. London Math.
  Soc. (3) \textbf{63} (1991), 165--211.

\bibitem{donkext}
S.~Donkin, \emph{On {${\mathrm{Ext}}^1$} for semisimple groups and
  infinitesimal subgroups}, Math. Proc. Camb. Phil. Soc. \textbf{92} (1982),
  231 -- 238.

\bibitem{donkquant}
\bysame, \emph{Standard homological properties for quantum {$\mathrm{GL}_n$}},
  J. Algebra \textbf{181} (1996), 235--266.

\bibitem{donkbk}
\bysame, \emph{The {$q$}--{S}chur {A}lgebra}, London Math. Soc. Lecture Note
  Ser., vol. 253, Cambridge University Press, Cambridge, 1998.

\bibitem{irv}
R.~S. Irving, \emph{The structure of certain highest weight modules for
  {${\mathrm{SL}}_3({K})$}}, J. Algebra \textbf{99} (1986), 438--457.

\bibitem{jandar}
J.~C. Jantzen, \emph{{D}arstellungen halbeinfacher {G}ruppen und ihrer
  {F}robenius--{K}erne}, J. {r}eine {a}ngew. Math. \textbf{317} (1980),
  157--199.

\bibitem{jantz2}
\bysame, \emph{{R}epresentations of {A}lgebraic {G}roups}, Mathematical surveys
  and monographs, vol. 107, AMS, 2003, second edition.

\bibitem{mathieu}
O.~Mathieu, \emph{Filtrations of {$G$}--modules}, Ann. Sci. \'Ecole Norm. Sup.
  (4) \textbf{23} (1990), no.~4, 625--644.

\bibitem{parker1}
A.~E. Parker, \emph{The global dimension of {S}chur algebras for
  {$\mathrm{GL}_2$} and {$\mathrm{GL}_3$}}, J. Algebra \textbf{241} (2001),
  340--378.

\bibitem{yehia}
S.~{el}~B. Yehia, \emph{{E}xtensions of simple modules for the universal
  {C}hevalley groups and its parabolic subgroups}, Ph.D. thesis, University of
  Warwick, 1982.

\end{thebibliography}

\providecommand{\bysame}{\leavevmode\hbox to3em{\hrulefill}\thinspace}
\providecommand{\MR}{\relax\ifhmode\unskip\space\fi MR }
\providecommand{\MRhref}[2]{%
  \href{http://www.ams.org/mathscinet-getitem?mr=#1}{#2}
}
\providecommand{\href}[2]{#2}

\end{document}